\input amstex

\documentstyle{amsppt}
\magnification=\magstep1
\nologo
\overfullrule=0pt
\hoffset=.25truein
\vsize=8.75truein
\mathsurround=2pt

\input epsf.sty

\topmatter
\title{Volume, diameter and the minimal mass of a stationary $1$-cycle}\endtitle
\author{A. Nabutovsky and R. Rotman}\endauthor
\rightheadtext{Nabutovsky and Rotman}
\leftheadtext{Nabutovsky and Rotman}
\endtopmatter
\address{Department of Mathematics, University of Toronto, Toronto,
Ontario, M5S 3G3, Canada.}\endaddress
\email{alex\@math.toronto.edu and rina\@math.toronto.edu}\endemail


\abstract
In this paper we present upper bounds on the minimal mass of a 
non-trivial stationary $1$-cycle.
The results that we obtain are valid for all
closed Riemannian manifolds.  The first result is that the minimal 
mass  of a stationary $1$-cycle on a closed $n$-dimensional 
Riemannian manifold $M^n$ is bounded from above
by $\frac{(n+2)!d}{3} $, where $d$ is the diameter of a 
manifold $M^n$. 
The second result is that the minimal mass of a stationary $1$-cycle 
on a closed Riemannian manifold $M^n$ is bounded
from above by 
$2 (n+2)! FillRad (M^n) \leq 2 (n+2)!(n+1)n^n\sqrt{n!} (vol(M^n))^{1/n}$,
where $FillRad (M^n)$ is the filling radius of a manifold, and
$vol(M^n)$ is its volume. 
\endabstract
\endtopmatter

\document
\baselineskip=16pt

{\bf 1. Introduction.}
\medskip
In 1983 M. Gromov asked whether there exists a constant $c(n)$
such that the length of a shortest closed geodesic 
$l(M^n) \leq c(n) (vol(M^n))^{1\over n}$, (see [G], p. 135).  
This problem also appeared as Problem 87 in a list of open problems
in Differential Geometry composed by S.-T. Yau ([Y], p. 689, or [SY], p. 297).
In the same spirit it might be interesting
to know whether  there exists $\tilde{c}(n)$ such that 
$l(M^n) \leq \tilde{c} (n) d$.  In this paper we prove the existence of
a stationary $1$-cycle such that its mass satisfies
these inequalities. 
In fact, our proofs demonstrate the existence of a stationary 1-cycle
of a special type which we will call {\it optimal}
such that its mass satisfies these inequalities. Optimal 1-cycles are defined
as follows: An optimal 1-cycle is a 1-cycle 
representable by a subset of $M^n$ that consists of finitely many
points $p_1,\ldots,p_l$ and a finite collection of
geodesic segments between some of these
points. Some of these segments can be (possibly trivial) geodesic loops.
We also allow
several different geodesic segments between two points, as well
as several copies of the same geodesic segment (in other words, every segment
can have a positive integer
multiplicity). Each point $p_i$ should be adjacent to at least
one of these segments. For any $p_i$
the number of the geodesic segments adjacent to $p_i$ counted with
their multiplicities is an even
positive number. (Here every geodesic loop
that starts and ends at $p_i$ is counted twice.)
This condition is necessary and sufficient in order for
the subset to represent a 1-cycle.  Equivalently, we can define
an optimal 1-cycle as a 1-cycle represented by a map of a finite
multigraph with vertices of positive 
even degree into $M^n$ such that each edge
is mapped into a (possibly non-minimizing) geodesic connecting
the images of the adjacent vertices.
\par
An optimal 1-cycle is {\it stationary}
if and only if the following property holds:
For any $i=1,\ldots l$ the sum of
unit vectors in $T_{p_i}M^n$ tangent to all geodesic segments
meeting at $p_i$ 
(counted with their multiplicities) is equal to zero.
Here each tangent vector is directed from $p_i$.
(The direction  from $p_i$ is formally defined as follows.
Parametrize the geodesic
segment. The direction  from $p_i$ coincides with the direction of the
first derivative, if the geodesic segment exits $p_i$, and is the opposite
direction,
if the geodesic segment enters $p_i$. It is easy to see that this definition
does not depend on the parametrization.)
The mass of an optimal 1-cycle is equal to
the sum of lengths of images of all its edges.
We also obtain an explicit upper bound for the total
number of all segments (counted with their multiplicities)
in an optimal stationary 1-cycle of mass not exceeding $\tilde c(n)diam(M^n)$
(or $c(n)vol(M^n)^{1\over n}$) in terms of $n$. 

Of course,
our estimates would give the estimates on the length of a shortest
closed geodesic if the stationary 1-cycle we
 obtain is realized by 
a closed geodesic. And, in fact,  when 
$M^n$ is diffeomorphic to the $2$-dimensional sphere,
this technique
produces a closed geodesic,
 as was observed by Pitts, Calabi and Cao,
(see [ClCo]). This fact enabled us to 
obtain estimates for
the length of a shortest closed geodesic on a manifold
diffeomorphic to $S^2$, (see [NR1]) improving previously
known results by C.B. Croke and M. Maeda, (see [C], [Ma]).
(Similar results were independently obtained by S. Sabourau in [S1].
Sabourau had also found curvature free upper bounds on the length of a 
shortest geodesic loop on a compact Riemannian manifold, (see [S2]).
Note that a geodesic loop is a stationary $1$-cycle if and only if it is 
a closed geodesic.
Also note that the general upper bounds for the length of shortest
closed geodesics can
be found in [NR]; see also the earlier paper [R]. Those estimates, however, 
use information about sectional curvature.)

Our techniques were inspired by the geometric  measure
theory approach to the existence of minimal submanifolds 
developed by Almgren and Pitts, (see [P]). We also use an appropriate
generalization of obstruction to an extension technique
used by Gromov in [G]
(see section 1.2 of [G] as well as the proof of Proposition
on p. 136 of [G]), and in the
case of the Theorem 2 below we use Gromov's upper bound
for the filling radius in terms of volume.

One of our starting points is 
the Almgren and Pitts minimax argument in the
geometric measure theory that produces stationary varifolds by considering
minima over non-trivial homotopy classes of maps $f$ of $S^k$ into the
spaces of cycles on a compact Riemannian manifold of $\max_{t\in S^k}$
mass$(f(t))$ (see Theorem 4.10 in [P]). (One gets the stationary varifold
for a $t\in S^k$ where the maximum of the mass of $g(t)$ is attained for a map 
$g$ minimizing  this maximum over the homotopy class.)
Roughly speaking one can get an upper bound for the volume of stationary
varifolds obtained by this method as follows:
One can start from an arbitrary map $f_0$ in a non-trivial
homotopy class. Denote by $M_0$ the maximal mass of cycles $f_0(t)$ over
$t\in S^k$. Then $M_0$ will provide an upper bound for the volume of
the varifold. 
\par
Non-trivial elements of homotopy groups of the space of cycles  correspond
to non-trivial homology classes of $M^n$ by virtue of
the classical result
of F. Almgren ([A]) that establishes the isomorphism between
the $k$-th homotopy group of the space of $l$-dimensional cycles $Z_l(M^n,G)$
in $M^n$ and
the $(k+l)$th homology group of $M^n$ with coefficients in $G$
(for all $k$, $l$ and $G=Z$ or $G=Z_p$ for some $p$).
The underlying geometric idea behind the
above correspondence is  slicing. For example,
consider a smooth function $f:M^n\longrightarrow [0,1]$ such that
$f^{-1}(\{0,1\})$ consists of finitely many points and such that 
$f^{-1}(t)$ is an $(n-1)$-dimensional cycle for any $t$. 
Then we obtain
a slicing of $M^n$ into $(n-1)$-dimensional cycles $f^{-1}(\{t\}), t\in [0,1]$
or, in other words, a map from $[0,1]$ to the space of $(n-1)$-dimensional
cycles in $M^n$. Since $0$ and $1$ are mapped to the zero cycle we
obtain, in fact, a map of $S^1=[0,1]/\{0,1\}$ into the space $Z_{n-1}(M^n)$.
It gives rise to an element of $\pi_1(Z_{n-1}(M^n))$ that corresponds to
the fundamental homology class of $M^n$ under the Almgren isomorphism,
and therefore is non-trivial. In fact, Almgren started from
assigning  a $(k+l)$-dimensional chain
(corr. a $(k+l-1)$-dimensional cycle) in $M^n$ with coefficients   
on $G$ to a map of $[0,1]^k$ (or $\partial [0,1]^k$)
into the space of cycles $Z_l(M^n,G)$. A simple modification of
this assignment which we will for brevity call {\it the Almgren
correspondence} will be used in this paper, and will be explained
in section 4 (in our context).
\par
In this paper we are interested only in one-dimensional cycles.
Therefore we are spared a significant amount of technical difficulties
arising in geometric measure theory. In particular, we can afford
the luxury of considering only 1-cycles
made of finitely many Lipschitz curves on $M^n$ (as in [ClCo]).
Therefore all of the above ideas
from the geometric measure theory 
become much more obvious and geometric
in our situation.

Now, let us denote the minimal mass (=length) of a non-trivial optimal
stationary
$1$-cycle in $M^n$
by $\alpha(M^n)$. Theorems 1 and 2 below establish 
upper bounds for $\alpha(M^n)$ in terms of the volume $vol(M^n)$
and the diameter $d$ of the manifold $M^n$.

\proclaim{Theorem 1}
Let $M^n$ be a closed simply-connected  Riemannian manifold
of dimension $n$. Let $q(\leq n)$ denote the minimal dimension $i$ such that
$\pi_i(M^n)\not =0$.

Then there  exists a non-trivial optimal stationary 1-cycle on $M^n$
that consists of
at most $\frac{(q+2)!}{2}$ geodesic segments and loops
such that its mass does not exceed
$\frac{(q+2)!}{3} d$.
If $q=2$ then there exists a stationary 1-cycle of length $\leq 4d$
that is either a closed geodesic or consists of two geodesic loops
emanating from the same point $p$. (Observe, that in the last case
the angles made
by the tangent vectors of these loops at $p$ are equal, and the bisectors of
these angles in $T_pM^n$ are oppositely directed rays.) 
\endproclaim

{\bf Remark.} It is well-known (and easy to prove) that if $M^n$ is {\it not}
simply-connected then there is a closed geodesic on $M^n$ of length $\leq 2d$.

In order to state the next theorem we will need the following definitions
and the following result of Gromov (see [G]).

\proclaim{Definition 1}
Let $M^n$ be a manifold topologicaly imbedded into an arbitrary metric space 
$X$. Then its filling radius, denoted $Fill Rad (M \subset X)$, is the
infimum of $\epsilon > 0$, such that $M^n$ bounds in the 
$\epsilon$-neighborhood $N_\epsilon(M^n)$, i.e. homomorphism 
$H_n(M^n,Z_2) \longrightarrow H_n(N_\epsilon(M^n),Z_2)$
induced by the inclusion map vanishes.
\endproclaim

\proclaim{Definition 2}
Let $M^n$ be an abstract manifold. Then its filling radius,
denoted $Fill Rad M^n$ will be $Fill Rad (M \subset X)$, where
$X=L^{\infty}(M)$, i.e. the Banach space of bounded Borel functions $f$
on $M^n$ and the imbedding of $M^n$ into $X$ is the map that 
to each point $p$ of $M^n$ assigns the distance function $p \longrightarrow
f_p =d(p,q)$
\endproclaim

\proclaim{Theorem A ([G])}
Let $M^n$ be a closed connected Riemannian manifold. Then
$Fill Rad M^n \leq (n+1) n^n n!^{\frac{1}{2}}(vol M^n)^{\frac{1}{n}}$.
\endproclaim

\proclaim{Theorem 2}
Let $M^n$ be a closed Riemannian manifold.  Then 
there exists a non-trivial optimal stationary 1-cycle
in $M^n$ made of at most $(n+2)!/2 $ geodesic segments and loops 
such that its mass is bounded from above by  
$(n+2)! Fill Rad M^n \leq (n+2)!(n+1)n!^{1\over 2}n^n(vol(M^n))^{\frac{1}{n}}$.
\endproclaim

{\bf 2. Basic definitions.} 
\medskip
First, let us briefly recall some basic notions of (the one-dimensional)
 geometric measure
theory. Our exposition was influenced by
the treatment of this subject in [ClCo], [BZ], [P] and [M].
We refer the reader to [M] or [P] for more details.

We say that a 1-dimensional
Borel subset $A$ of $M^n$ is {\it rectifiable} if
there exists a Borel subset $S\subset R^1$ of finite measure and a $1-1$
Lipschitz map $f:A\longrightarrow M^n$ such that $f(S)\subset A$ and the
1-dimensional Hausdorff measure of $A\setminus f(S)$ is zero.
Any such set $A$ can be regarded as a linear functional
on the space of $C^\infty$-smooth 1-forms on $M^n$: $T_A(\phi)=\int_A\phi$.
When $A$ is regarded  as the linear functional $T_A$ on 1-forms,  
it is   called  a {\it rectifiable current}.
In general, 1-dimensional
rectifiable currents are defined as continuous linear functionals on
the space of 1-forms on $M^n$ (with the mass topology explained below)
that can be represented as countable linear combinations $\Sigma_{i=1}^{\infty}
a_iT_{A_i}$ with integer
coefficients $a_i$ of functionals $T_{A_i}$ corresponding to 1-dimensional
rectifiable sets $A_i$. It is not difficult to see that an arbitrary
rectifiable 1-current can be represented as $\Sigma_{i=1}^{\infty}a_iT_{M_i}$,
where $a_i$ are integer coefficients, and $M_i$ are 1-dimensional
submanifolds of $M^n$, such that
$\Sigma_{i=1}^{\infty} \vert a_i\vert l(M_i)<\infty$,
where $l$ denotes the length.

The space of rectifiable 1-currents can be endowed with the weak topology:
$T_i\longrightarrow T$
if and only if $T_i(\phi)\longrightarrow T(\phi)$ for any 1-form $\phi$.
This space can also be endowed with the the mass norm as follows:
Define the {\it mass norm} of a 1-form $\omega$ as
$\sup_{x\in M^n}\sup_{v\in T_xM^n; \Vert v\Vert=1}\vert \omega(x)(v)\vert$.
The {\it mass} $M(T)$ of a 1-dimensional rectifiable current $T$ is defined
as the supremum over the set of all 1-forms $\omega$ on $M^n$ of mass one
of $\vert T(\omega)\vert$. Note that if $A$ is a $C^1$-smooth curve in $M^n$
with a countable set of self-intersections then $M(T_A)=l(A)$. However,
if $A$ backtracks over itself than the lengths of its pieces
where the bactracking takes place cancel each other, and the mass is
strictly less than the length. The mass is a lower-semicontinuous
(but not continuous) functional on the space of rectifiable currents with
the weak topology. The topology induced by the mass norm on the space of
1-currents is called the {\it mass topology}. Observe that a $C^1$-smooth
path in the space of smooth closed curves in $M^n$ regarded as currents is
always continuous in the weak topology but usually NOT continuous 
in the mass topology. Indeed, the mass of the difference of two $C^1$-close
but not intersecting $C^1$-smooth closed curves is equal to the {\it sum}
of their lengths. (More generally, if $A_i$ are Lipschitz curves in $M^n$
such that the sets of their pairwise intersections and self-intersections
have measure zero, then $M(\Sigma_i a_iT_{A_i})=\Sigma_i \vert a_i\vert
l(A_i)$.)
A sequence of 1-currents corresponding to Lipschitz curves $\gamma_i$
converges to the 1-current corresponding to a Lipschitz curve
$\gamma$ in the mass topology if and only if the mass of the
symmetric difference
$(\gamma\setminus\gamma_i)\bigcup (\gamma_i\setminus\gamma)$
converges to zero. 

Let $T$ be a 1-dimensional
rectifiable current. If for any smooth function $f$ on $M^n$
$T(df)=0$, then we say that $T$ is a {\it rectifiable 1-cycle}
on $M^n$. Define $\bar Z_1(M^n,Z)$ as the subspace of the space of all
1-dimensional rectifiable currents endowed with the weak topology that consists
of all 1-cycles. Consider the space of 1-cycles made of finitely many
curves $Z_1(M^n,Z)\subset \bar Z(M^n,Z)$ that consists
of all 1-cycles made of finitely many Lipschitz curves (in other words, we are
considering only those cycles that can be represented as $\Sigma_{i=1}^k
a_iT_{A_i}$, where $A_1,\ldots A_k$ are images of $[0,1]$ under Lipschitz
maps to $M^n$,
 and $a_i\in Z$.) Any such cycle is an integer linear combination of currents
corresponding to {\it closed} Lipschitz curves in $M^n$ with integer
multiplicities. (Note that the above notation is not standard:
the notation $Z_1(M^n,Z)$ is usually used for the space that we denote
$\bar Z_1(M^n,Z)$ in this paper.) The approximation theorem from the
geometric measure theory (cf. [M]) implies that $\bar Z_1(M^n,Z)$ is the
closure of $Z_1(M^n,Z)$ in either weak or mass topology.
\par
Now we would like to introduce some spaces of ``nice" 1-cycles
that are especially useful for our purposes.
Following [ClCo] it is convenient to consider spaces of {\it parametrized}
1-cycles made of at most  $k$ closed curves: Define $\Gamma_k$ as
the space of all $k$-tuples $(\gamma_1,\ldots ,\gamma_k)$ of Lipschitz
maps of $[0,1]$ to $M^n$ such that $\Sigma_{i=1}^k\gamma_i(0)=\Sigma_{i=1}^k
\gamma_i(1)$. Endow $\Gamma_k$ with the following metric topology:
First, because of the Nash embedding theorem we can assume without any loss
of generality that $M^n$ is isometrically
embedded into the Euclidean space $R^N$ of a large
dimension. Now define the distance
$d((\alpha_1,\ldots,\alpha_k),(\gamma_1,\ldots,\gamma_k))$ as
$\max_{i,t}d_{M^n}(\alpha_i(t),\gamma_i(t))+\Sigma_{i=1}^k\sqrt{\int_0^1 \vert
\alpha'_i(t)-\gamma'_i(t)\vert^2 dt}$.
It is easy to see that the length functional $l((\gamma_1,\ldots,
\gamma_k))=\Sigma_{i=1}^k l(\gamma_i)$ is a continuous functional on this
space. Also, the map $I:\Gamma_k\longrightarrow Z_1(M^n,Z)$ (with the 
weak topology) given
by the formula $I((\gamma_1,\ldots,\gamma_k))=\Sigma_{i=1}^k T_{\gamma_i}$
is continuous. Denote the image of $\Gamma_k$ under $I$ by
$Z_{(k)}$. We will call $Z_{(k)}$ {\it the space of 1-cycles on $M^n$ made
of at most $k$ closed curves}, and we will call $k$ {\it the order} of
of these 1-cycles. 
Observe that $\Gamma_k$ contains all collections of
at most $k$ suitably parametrized closed curves in $M^n$.
Therefore $Z_{(1)}\subset Z_{(2)}\subset Z_{(3)}\subset\ldots$, and
$\bigcup_{i=1}^\infty Z_{(i)}=Z_1(M^n,Z)$.
Also, for any $x$ let the subset of $\Gamma_k$ formed by all $\gamma=
(\gamma_1,\ldots ,\gamma_k)$ such that $l(\gamma)=\Sigma_{i=1}^k l(\gamma_i)
\leq x$ be denoted by $\Gamma_k^x$, and   the image of $\Gamma_k^x$
under $I$ be denoted by $Z_{(k)}^x$.
We will call $Z_{(k)}^x$ {\it the space of 1-cycles
on $M^n$ of length $\leq x$ made of at most $k$ curves}.
Similarly, we will call elements of $\Gamma_k$ {\it parametrized 1-cycles
made of $k$ curves}, $k$ will be called
{\it the order} of parametrized cycles from $\Gamma_k$,
and elements of $\Gamma_k^x$ will be
called {\it parametrized 1-cycles of length $\leq x$
made of $k$ curves}.
\par
Our definition of {\it stationary} parametrized or non-parametrized
1-cycles made of $k$ segments is quite similar to the standard
definition in the geometric measure theory (cf. [P]):
Let $X$ be a vector field on $M^n$. It determines the one-parameter
group of diffeomorphisms $\Phi_X(t)$ of $M^n$.
For any $\gamma\in\Gamma_k$
$\Phi_X(t)(\gamma)$ is a continuous one-parameter family of parametrized
1-cycles. Now we can consider the continuos function $L_{X,\gamma}(t)$
defined as the total length of $k$ Lipschitz curves
that together form $\Phi_X(t)(\gamma)$. If for any $X$ $t=0$ is the
critical point of $L_{X,\gamma}$ then $\gamma$ is called {\it a stationary
parametrized 1-cycle made of $k$ geodesic segments}. Its image
$I(\gamma)$ in $Z_{(k)}$ is called {\it a stationary (non-parametrized)
1-cycle made of $k$ geodesic segments}. The value of
the derivative of $L_{X,\gamma}$ at $t=0$ will be called {\it the first
variation of the length of $\gamma$ in the direction of $X$}. 
In order to understand the stationarity condition we first
observe that each of $k$ curves $\gamma_i$ comprizing $\gamma$ must
be a geodesic in a neighborhood of every point $\gamma_i(t)$, where
$\gamma_i$ does not intersect itself or $\gamma_j$ for some $j\not= i$,
if $\gamma$ is a stationary parametrized 1-cycle.
(Otherwise we can use a local vector field $X$ supported on a small
ball centered at this point of $\gamma_i$ to demonstrate
that $\gamma$ is not stationary.) Now we can consider vector fields
supported in a small neighborhood of $a=\gamma_i(1)$ or $a=\gamma_i(0)$.
When the neighborhood is small the first variation of length is dominated
by $<X(a),\Sigma_{i\in I_1}\gamma'_i(1)-\Sigma_{i\in I_2}\gamma'_i(0)>$,
where $I_1$ is the set of $i\in\{1,\ldots,k\}$ such that $\gamma_i(1)=a$,
and $I_2\in\{1,\ldots,k\}$ is the set of $i$ such that $\gamma_i(0)=a$.
Hence for any such $a$ the sum of the tangent vectors to curves $\gamma_i$
meeting at $a$ (and directed to $a$)
must be equal to zero.
A similar condition must hold at any point $\gamma_i(t)$, where $\gamma_i$
has an intersection of finite multiplicity with itself and/or other
curves $\gamma_j,\ j\not =i$. However, in general,
 $\gamma_i$ need not be a geodesic in a neighborhood of such a point
even if $\gamma$ is stationary. For example, let $k=1$. Assume that
$\gamma=\gamma_1$ is a three-petal curve that
consists of three geodesic loops emanating from the same point $p$.
Consider three angles formed by tangent vectors at $p$ to these three loops.
(As usual in this paper, we direct these tangent vectors from $p$.)
Assume further
that: 1) These three angles have equal values that are strictly
less than $\pi/6$;
2) The bisectors of these three angles lie in a plane in
$T_pM^n$ and form angles equal to $\pi/3$ with each other. Then it is easy
to see that $\gamma\in \Gamma_1$
is a stationary parametrized 1-cycle, but is not
a closed geodesic (and does not correspond to a closed geodesic in any obvious
way). This example motivates the following definition: We say that
a parametrized 1-cycle $\gamma\in \Gamma_k$
is {\it strongly stationary} if it is
stationary, and for any $i\in\{1,\ldots ,k\}$ $\gamma_i$ is a geodesic.
If a 1-cycle $z$ is equal to $I(\gamma)$ for some $k$ and some strongly
stationary parametrized 1-cycle $\gamma\in\Gamma_k$, then we say that
$z$ is a {\it strongly stationary} non-parametrized
1-cycle. From now on we will be using only the notion of strong stationarity.

It is easy to see that our
definition of a non-parametrized strongly
stationary 1-cycle is equivalent
to the defintion of a stationary optimal 1-cycle
given in the introduction. Also note that a strongly stationary
1-cycle made of one geodesic segment must be a geodesic loop and therefore
must be a closed geodesic. (Two unit vectors tangent  to the geodesic loop
at its origin must cancel each other.) A strongly
stationary 1-cycle made of two
geodesic segments either consists of two geodesic segments connecting
two different points or consists of two geodesic loops. In the first case
it is easy to see that it is a closed geodesic. In the second case
we have two subcases. If these geodesic loops are based at 
different points, then both of them must be closed geodesics.
If they are based at the same origin, then
we see that the sum of the four unit tangent vectors at the
origin
of the loops must be equal to zero.
If the dimension of the manifold
is equal to $2$, then this condition implies that this strongly stationary
1-cycle
is just a self-intersecting closed geodesic.
If the dimension of
the manifold is greater than two, then, in principle, these two geodesic
loops need not form a closed geodesic. But it is easy to see that 1) the
angles between the tangent vectors to these loops at their common origin
$p$ are equal; 2) The bisectors of these two angles in $T_pM^n$ are
oppositely directed rays. This condition is, in fact, equivalent to
the strong
stationarity for a 1-cycle formed by two geodesic loops emanating from
the same point.
\medskip
{\bf 3. A Morse-theoretic type lemma for $\Gamma_k$}
\medskip
The main technical results of this section
resemble Theorem 4.3 in [P] (though they do not directly
follow from it).
On the other hand, they resemble a basic
result from the Morse theory asserting that if there are no critical
points of a smooth function $F:M\longrightarrow R$ on
a compact manifold $M$ in the set $F^{-1}([x_1,x_2])$ then the sublevel set
$F^{-1}((-\infty, x_1])$ is a deformation retract of the sublevel set
$F^{-1}((-\infty, x_2])$. (The deformation retraction can be obtained
using the gradient flow of $F$.) Our goal is to obtain a result of such type
for the length functional on $\Gamma_k$. The main technical problem is that
$\Gamma_k$ is not an infinite-dimensional manifold, but consists of finitely
many intersecting pieces (each of which is an infinite-dimensional manifold).

\proclaim {Lemma 3}
Assume that there are no non-trivial strongly
stationary 1-cycles on $M^n$ of length
$\leq x$  made of $k$ geodesic segments. 
Let $f:S^i\longrightarrow \Gamma_k^x$ be a continuous map.
Then
there exists a homotopy $H:S^i\times [0,1]\longrightarrow \Gamma_k^x,
H(x,0)=f(x), H(x,1)\in \Gamma_k^0\subset \Gamma_k^x$ between
$f$ and a map $g$ of $S^i$ into $\Gamma_k$
such that for any $x\in S^i$
$g(x)$ consists of $k$ constant maps of $[0,1]$ into $M^n$.
\endproclaim

Observe that the conclusion of Lemma 3 falls short of the desired
contractibility of $f$. However, if we project everything to the space
of non-parametrized 1-cycles, then the contractibility is immediate.
Also, we can prove the contractibilty if some auxiliary conditions
on $f$ are satisfied. More precisely, we have the following proposition:

\proclaim {Proposition 4} 
Assume that there are no non-trivial strongly
stationary 1-cycles on $M^n$ of length
$\leq x$ made of $k$ geodesic segments.
Let $f:S^i\longrightarrow Z^x_{(k)}$ be a continuous map such that
there exists a map homotopic
to $f$ that can be lifted to a map $F:S^i\longrightarrow \Gamma_k^x$. (That is,
$f$ is homotopic to the composition of $I$ and $F$.) Then
$f$ is contractible (in $Z^x_{(k)}$). Further,
consider the map $F_j:S^i\longrightarrow M^n$
for every $j\in\{1,\ldots,k\}$
defined by the formula $F_j(p)=(F(p))_j(0.5)$ for any $p\in S^i$.
If all these maps $F_j$ are contractible,
then $F$ is contractible.
\endproclaim

{\bf Proof of Proposition 4 assuming Lemma 3:}
Without any loss of generality we can assume that $f$ is the composition of 
$I$ and $F$.
Apply Lemma 3 to $F$.
There exists a homotopy $H$ between $F$ and a map of $S^i$ into
$\Gamma_k^0$.
Observe that $I$ sends all elements of $\Gamma_k^0$ to the zero 1-cycle.
Therefore the composition of $I$ and $H$
constitutes a homotopy between $f$ and the constant
map of $S^i$ into the zero 1-cycle, i.e. a contraction of $f$.
\par
It remains to verify that if $F_j$ is contractible for each $j$, then
$H(1):S^i\longrightarrow \Gamma_k^0$ is contractible. Note that
$H(1)$ can be regarded as a map of $S^i$ into $(M^n)^k$.
Therefore it is sufficient to check only that $k$ maps $H(1)_j$
are contractible. But each of these maps can be connected with
$F_j$ by a homotopy $Q$ that can be defined 
by the formula $Q(t)(p)=(H(t)(p))_j(0.5)$ for any $t\in [0,1]$,
$p\in S^i$.
QED.

{\bf Proof of Lemma 3:}

For the convenience of
the reader the proof of Lemms 3 will be split into several steps.
\par
{\bf 3.1.}
We are going to demonstrate the existence of a deformation of $\Gamma_k^x$
into $\Gamma_k^0$, which clearly implies the lemma.
Recall that a deformation of a space $X$ into its subset $A$ is, 
by definition, a continuous map $g$ of $X$ into itself homotopic to
the identity map $X\longrightarrow X$ such that $g(X)\subset A$. More
generally,
if $B$ is a subset of $X$ such that $A\subset B$ we say that a map
$g:B\longrightarrow A$ is a deformation in $X$ if there exists a homotopy
$G:B\times [0,1]\longrightarrow X$ such that for any $b\in B$ $G(b,0)=b$ and
$G(b,1)=g(b)$.
But sometimes by a deformation
we will mean the whole homotopy $G$ and not just $g(b)=G(b,1)$. 
\par {\bf 3.2. Birkhoff curve-shortening process for $\Gamma_k^x$.}
First, we are going to proceed as in  the first step of the Birkhoff
curve shortening process described in [C] or [ClCo]: Let $inj(M^n)$
denote the injectivity radius of $M^n$. Choose
$N=[4x/inj(M^n)]+ 1$.
Let $\gamma$ be an element of $\Gamma_k^x$.
Divide each of $k$ segments $\gamma_i$ of $\gamma$
into $N$ pieces of equal length by points $\gamma_i(t_{ij}),
j\in\{0,1,\ldots ,N\}, t_{i0}=0, t_{iN}=1$.
Consider the unique minimizing
geodesic segments between 
$\gamma_i(t_{ij})$ and $\gamma(t_{i,j+1})$ for all $j$.
The length of each of these segments  does not exceed $inj(M^n)/4$.
For any $i$ $N$ such geodesic segments form a piecewise geodesic $\bar\gamma_i$
connecting $\bar\gamma_i(0)=\gamma_i(0)$ with  $\bar\gamma_i(1)=\gamma_i(1)$.
The length of $\bar\gamma_i$ does not exceed the length of $\gamma_i$. There
exists the following homotopy between $\gamma_i$ and $\bar\gamma_i$:
At the moment of time $\tau\in [ 0,1]$ we follow each of the $N$ segments
of $\gamma_i$ (between $\gamma_i(t_{ij})$ and $\gamma_i(t_{i,j+1})$) from
$\gamma_i(t_{ij})$ to $\gamma_i(\tau t_{ij} + (1-\tau) t_{i,j+1})$ and then
make a shortcut from $\gamma_i(\tau t_{ij} + (1-\tau) t_{i,j+1})$ to
$\gamma_i(t_{i,j+1})$ along the shortest geodesic. Then we reparametrize
the resulting curve proportionally to its arc length.
It is easy to see that the length of the curve does not increase during this
homotopy. 
Denote the resulting homotopy by $h_i$.
Combining all
of these $k$ homotopies $h_i$
we obtain a homotopy between the parametrized 1-cycle
$\gamma$ and the parametrized 1-cycle $\bar\gamma=(\bar\gamma_1,\ldots ,
\bar\gamma_k)$ made of $k$ piecewise geodesics with $N$ breaks.
This homotopy
depends on $\gamma$ in a continuous way. Therefore we obtain a deformation
of $\Gamma_k^x$ into its subset $g_{k,N}^x$ defined as the set of
all parametrized 1-cycles $\bar\gamma=(\bar\gamma_1,\ldots ,\bar\gamma_k)$
from $\Gamma_k^x$ such that for any $i$ $\bar\gamma_i$ is 
a piecewise geodesic made of $N$
geodesic segments of non-zero length $\leq inj(M^n)/4$
parametrized proportionally to its arclength.
Let us denote this deformation by $B_N$.
Below we will refer to it as the {\it the Birkhoff deformation}.
We regard $g_{k,N}^x$ as the
subset of a larger set $G_{k,N}^x$ defined as the set of 
all elements $\gamma=(\gamma_1,\ldots ,\gamma_k)$ of $\Gamma_k^x$ such 
that for any $i$ $\gamma_i$ is a piecewise geodesic made of at most $N$
geodesic segments of non-zero length $\leq inj(M^n)/2$. In other words,
the only difference between $g_{k,N}^x$ and $G_{k,N}^x$ is that we allow
elements of $G_{k,N}^x$ to have
somewhat longer geodesic segments.
\par 
Now we are going
to prove that $g_{k,N}^x$ can be
deformed into its subset $g_{k,N}^0=\Gamma_k^0$ inside $\Gamma_k^x$.
\par {\bf 3.3 A classification of elements of $\Gamma_k^x$ and $G_{k,x}^N$.}
Let us define an equivalence relation on $\Gamma_k^x$.
For any element $\gamma=(\gamma_1,\ldots,\gamma_k)$ of
$\Gamma_k^x$ or $G_{k,N}^x$ consider
its $2k$ endpoints $\gamma_i(0),\gamma_i(1)$.
We will call these points {\it multiple points} of $\gamma$.
The set of these
$2k$ points can be partitioned into $J$ non-empty sets $A_j$,
($J\in\{1,\ldots,k\}$), such that 1) Each set $A_j$ contains the
equal number of points of the form $\gamma_i(0)$ (for some $i$) and
$\gamma_l(1)$ (for some $l$); 2) $\gamma_i(t_1)=\gamma_l(t_2)$
for some $i,l\in\{1,\ldots,k\}$ and $t_1,t_2\in\{0,1\}$ if and only if
$\gamma_i(t_i)$ and $\gamma_l(t_2)$ are in the same set $A_j$ for some
$j\in\{1,\ldots,J\}$.
The number $J$ will be
called {\it the number of multiple points} of $\gamma$.
We will say that two 1-cycles $\gamma$ and $\beta$
from $\Gamma_k^x$ are of the same {\it type} if these partitions 
for $\gamma$ and $\beta$ coincide,
and the set of all $i$ such that $\gamma_i$ is constant coincides
with the set of values of $i$ such that $\beta_i$ is constant.
\par
For example, let $k=2$. In these
case there will be three types of parametrized $1$-cycles from $\Gamma_k^x$
when neither $\gamma_1$ nor $\gamma_2$ is constant:
(a) $\gamma_1(0)=\gamma_1(1)\not =\gamma_2(0)=\gamma_2(1)$ ($1$-cycle that
consists of 2 closed curves that do not intersect at their endpoints);
(b) $\gamma_1(0)=\gamma_2(1)\not =\gamma_2(0)=\gamma_1(1)$ ($1$-cycles
that consist of one closed curve obtained by glueing together $\gamma_1$
and $\gamma_2$, where endpoints of $\gamma_1$ (and of $\gamma_2$) are
different; and (c) $\gamma_1(0)=\gamma_1(1)=\gamma_2(0)=\gamma_2(1)$
($\gamma_1$ and $\gamma_2$ are loops with the common endpoint. Such $1$-cycle
can be considered as either made of two closed curves or of one
closed curve with the self-intersection.) There will be two types when both
$\gamma_1$ and $\gamma_2$ are constant, namely, $\gamma_1\not= \gamma_2$
and $\gamma_1=\gamma_2$. For each $i=1,2$ there will
also be two types corresponding to the case when $\gamma_i$ is constant
and $\gamma_{3-i}$ is a non-constant loop: One type corresponds
to the case when $\gamma_i(t)=\gamma_{3-i}(0)=\gamma_{3-i}(1)$ and the other
type corresponds to the case $\gamma_i(t)\not= \gamma_{3-i}(0)=\gamma_{3-i}(1)$.
\par
However, we will also
need a stronger equivalence relation on $G_{k,N}^x\subset
\Gamma_k^x$.
For each $\gamma\in G_{k,N}^x$
we will consider
$kN$
geodesic segments $\gamma_{ij}$
in $k$ parametrized curves $\gamma_i$ forming $\gamma$. Consider $(N-1)k$
endpoints of these geodesic segments  that are
not endpoints of $k$ curves $\gamma_i$.
We will call them {\it double points}
in order to emphasize that there are exactly two geodesic segments
meeting at any of these points. (However note, that it is possible that
precisely two geodesic segments meet at some of the multiple points as well).
We will say that two elements $\alpha,\beta$
of $G_{k,N}^x$ are of the same {\it type}
as elements of $G_{k,N}^x$ if they are of the same type as elements
of $\Gamma_k^x$ and for each
$i=1,\ldots ,k$, $j=1,\ldots ,N$ the geodesic segment $\alpha_{ij}$ is constant
if and only if $\beta_{ij}$ is constant.
Also note that the type of $\gamma$ as an
element of $G_{k,N}^x$ and the (vectors of) positions of $J$ multiple points
and $(N-1)k$ double points
on the manifold
determine an element of $G_{k,x}^N$ uniquely.
Therefore for any specific type 
of $\gamma$ we can identify an element $\gamma=(\gamma_i)_{i=1}^k$ of
$G_{k,N}^x$ with the corresponding
element of $(M^n)^{J+(N-1)k}$,
where $J$ is the number of multiple points of $\gamma$.
Note that copies of
$(M^n)^{J+(N-1)k}$ corresponding to individual types
can be simultaneously embedded into the ambient
manifold $(M^n)^{(k+1)N}$ by the corresponding diagonal embeddings
that identify endpoints of $k$ intervals in accordance with the type.
\par {\bf 3.3.A A partial order on types of elements of
$G_{k,N}^x$.} We will say that a type $A$ is {\it higher}
than a type $B$, if: 1) Multiple points of $A$ can be obtained by merging
some of the multiple points of $B$. (In other words
the partition considered in 3.3 for $A$ is coarser than the partition for
$B$.); 2) If for some 1-cycle $\beta$ of type $B$, for
$i=1,2,\ldots ,k$, and for some
$j=1,\ldots N$  the geodesic segment $\beta_{ij}$
is constant,
then for any element $\gamma$ of type $A$ $\gamma_{ij}$ is also constant.
\par
The resulting relation is a partial order on the set of types.
For example, the maximal types
among types that correspond to cycles of non-zero length 
are the types with just one multiple point $m$,
$Nk-2$ constant geodesic segments,
and just two non-constant geodesic segments. These two segments
correspond to the minimal geodesic connecting $m$ with the only 
double point, $d$, and traverse this geodesic in opposite directions.
\par {\bf 3.4. Some general remarks about deformations of $g_{k,N}^x$
that will be constructed below.}
These deformations will consist of finitely many steps. Each
of these steps constitute either 
the Birkhoff deformation described in 3.2, or
will be a deformation inside $G_{k,N}^x$. In the last case in order
to describe the deformation we need to describe the trajectories of
individual multiple and double points.
More formally, during a deformation of
the last type we will be deforming each individual cycle $\gamma$
as an element of $(M^n)^{(N+1)k}$.
Each multiple or double point will move along trajectories
obtained as a projection of trajectories
of a vector field defined in an open set in $(M^n)^{(N+1)k}$
that contains the union of all
diagonnaly embedded copies of $(M^n)^{J+(N-1)k}$ corresponding
to various types of elements of $G_{k,N}^x$.
\par
We will construct this vector field in such a way that will not allow
the type of the elements of $G_{k,N}^x$ to change during these stages of
the deformation, but
at the very end of the deformation, when the trajectory of
the flow hits $G_{k,N}^0$, and the length of the element becomes
zero.
This is achieved by defining the components of the vector fields
that correspond to the individual multiple
and double points of $\gamma$ so that in the situation when the distance
between two distinct multiple and/or double points 
becomes small,
they will move along the trajectories of the same smooth vector field
on $M^n$. (Here we are talking only about pairs of points 
the collision of which can change the type.)
Moreover, since the type of any 1-cycle from $g_{k,N}^x$ regarded as
an element from $\Gamma_k^x$ cannot change during the Birkhoff deformation
(since all multiple points remain unchanged), $\Gamma_k^x$ types
of 1-cycles remain unchanged through the whole deformation
until its very last moment.
We will not need this feature of our construction in the proof 
of Lemma 3, but it will turn out to be convenient for the next section.
\par
However, we can encounter the following problem:
We need to ensure that the distance between
any two (double or multiple)
points that should be connected by a geodesic segment is less than
$inj(M^n)$. We will resolve this technical complication
in the following way.
Our choice of $N$ ensures that at the beginning these distances
do not exceed $inj(M^n)/4$ for any $\gamma\in g_{k,N}^x$.
Therefore we can deform our parametrized
1-cycles using a flow that will be constructed below
for a certain safe amount of time which is
sufficiently small to ensure that lengths of segments grow by
not more than $inj(M^n)/4$.
Then we
stop and apply
the Birkhoff deformation with $N$ breaks $B_N$ again.
That is, we forget that the $k$ curves forming an element of $G_{k,N}^x$
are already piecewise geodesics
with $N$ breaks, divide them into $N$ equal pieces (of length $\leq x/N\leq
inj(M^n)/4$) by $N-1$ points 
and replace the curve by another piecewise geodesic consisting
of $N$ geodesic
segments of length $\leq inj(M^n)/4$ with vertices in these
points. (Of course, this replacement is made using the
length non-increasing
homotopy described in 3.2 that was used in the definition of
the Birkhoff deformation.)
At this point, in principle, the type of $\gamma$
as an element of $G_{k,N}^x$ can change because new double points
can merge in a different pattern from what was before the Birkhoff stage.
(Note however that the $\Gamma_k^x$ type remain unchanged since $B_N$
does not affect multiple points.)
Now we are ready to continue the
deformation using the same flow again, etc... 
The vector fields and the times of these deformations will be chosen so that
at each stage of the
deformation every element of $G_{k,N}^x$ will become shorter
by at least
a certain $\delta=\delta(M^n,N,k,x)>0$. Therefore we will need only a finite
number of stages to reach $G_{k,N}^0$.
\par {\bf 3.5. The direcion of the steepest descent.}
We will start the construction of the flow from the following observation:
Since there is no strongly stationary
1-cycle of positive length $\leq x$, then for any piecewise
geodesic 1-cycle $\gamma$ from $G_{k,N}^x$
there exists a system of vectors at all
multiple and double points such that a small
deformation of $\gamma$ (regarded as an element of
$(M^n)^{J+(N-1)k}$) in
the direction of these vectors leads to an element of $G_{k,N}^x$ of
the same type but of a smaller length. These vectors are constructed
as follows:
Any multiple point corresponds to a set $A_j$ of the partition. 
The vector at this point is calculated as $\Sigma_{\gamma_l(t_i)\in A_j;
t_i=0\ or\ 1}v_l(t_i)$, where $v_l(t_i)$ is the unit vector tangent to
$\gamma_l$ at $t_i$ directed from the multiple point.
If we will regard all $Nk$ geodesic segments as curves in $M^n$, then
each double point 
is adjacent to two geodesic segments of the curve (with the exception
of the case, when this double point is connected by a sequence of geodesic
segments of zero length with a multiple point. In this last case
the double point coincides with the multiple point, and we {\it define}
the component of $v$ corresponding to this double point to be equal
to the already determined component of $v$ corresponding to
the multiple point.)
The component of $v$
at this double point is calculated as the sum of two unit vectors
tangent to two
geodesic segments that meet at this point and are directed from it.
Our assertion
now follows directly from the first variation formula 
for the length functional.
We will call the system of 
$J+(N-1)k$ tangent vectors of $M^n$ {\it a
deformation vector} for $\gamma$ and will denote it by $v(\gamma)$. Note that
$v(\gamma)$ is the collection of zero vectors if and only if $\gamma$
is a strongly stationary 1-cycle.
\par
The first variation of the length of $\gamma$
in the direction of $v(\gamma)$ is equal to $-\Vert v(\gamma)\Vert^2$,
where $\Vert v(\gamma)\Vert^2$ is 
calculated as follows:
When we deal with $J$ components of $v$ corresponding to multiple
points we just sum their squares.
We are going to say that a pair of double points {\it merges}
if they are connected by a sequence of {\it constant} geodesic segments.
We define a {\it cluster} of double points as a maximal set of double
points such that each pair of them merges.
We say that a cluster is {\it negligible} if one of double points in the
cluster is connected by a sequence of {\it constant} segments with
a multiple point (so geometrically, all double points in this cluster
coincide with the multiple point).
Our definition of $v$ implies
that components of $v$ corresponding to double points
in a cluster point are equal. When we calculate $\Vert v(\gamma)\Vert^2$
we {\it by definition} disregard all double points in each negligible cluster
and count the squared norm of the component of $v$ corresponding to
all double points in a cluster only once for each non-negligible cluster.
In other words we
define $\Vert v(\gamma)\Vert^2$ as $\Sigma_{m_i; i=1,\ldots J}
\Vert v(\gamma)(m_i)\Vert^2+\Sigma_{Non-negligible\ clusters\ of\ double\ 
points}\Vert v(\gamma)(d_i)\Vert^2$, where the first summation is over the set
of all
multiple points and the second summation is over the set of all non-negligible
clusters of double points.
\par
{\bf 3.6. The deformation vector for $\gamma$ can be used
to decrease the length for all
$\gamma_*$ sufficiently close to $\gamma$.}
Unfortunately, the dependence of $v(\gamma)$ on $\gamma$
is not continuous. This happens because the type of element of $G_{k,N}^x$
changes in a discontinuous manner. 
Yet, it is easy to see that for
any $\gamma\in G_{k,N}^x$ there exists a sufficiently small 
 open neighborhood $U$ of $\gamma$ in
$G_{k,N}^x$ and a
positive $\mu$ such that for any $\gamma_*\in U$ a sufficiently small
deformation
of $\gamma_*$ in the direction of the deformation 
vector $T_\gamma^{\gamma_*}(v(\gamma))$ defined below
decreases the length
of $\gamma_*$, and the first variation of the length in the
direction of $T_\gamma^{\gamma_*}(v(\gamma))$
is less than or equal to
$-\Vert v(\gamma)\Vert^2/2$.
The above deformation vector, $T_{\gamma}^{\gamma_*}(v(\gamma))$ is 
obtained from $v(\gamma)$
by the parallel transport of all $J+\Sigma_i (N_i-1)$ vector components
of $v(\gamma)$ along the shortest geodesics connecting 
the vertices of $\gamma$ with the corresponding vertices
of $\gamma_*$. 

More precisely, one first chooses $U$ so small that:
1) different multiple or double points of $\gamma$ cannot merge in $U$;
2) For any $\gamma_*\in U$ any of its multiple points has the unique closest
multiple point of $\gamma$, and any double point of
$\gamma_{*i}$ has the unique closest double or multiple
point of $\gamma_i$ at the distance not exceeding $inj(M^n)/4$.
But note that, in principle,  each multiple point of
$\gamma$ can split into two multiple
points (that can be connected or not connected by a geodesic segment)
or into a pair multiple point - double point for an arbitrarily small
$U$. Also note that if $N_i <N$ then each double
point of $\gamma_i$ can bifurcate into two distinct double points connected
by a very short geodesic.
If one of $k$ segments $\gamma_i$ is a constant geodesic loop,
then the corresponding
multiple point can bifurcate into a pair of points that consists of the 
multiple point and a double
point that is very close to the multiple point.
These points 
 are connected by two oppositely oriented copies of the 
shortest geodesic, together
forming a short piece-wise geodesic loop based at the multiple point.
Finally note that a finite number of bifurcations of these types
can occur simultaneously. So, the dimension of $T_\gamma^{\gamma_*}(v(\gamma))$
can be greater than the dimension of $v(\gamma)$.
But condition 2) in the definition of
$U$ implies that even if a multiple or a double point of $\gamma$
bifurcates into a finite number of (multiple and/or double)
points we know unambigiously 
how to define the corresponding component of
$T_\gamma^{\gamma_*}(v(\gamma))$ for each of them: we just perform 
the parallel transport
of the corresponding component of $v(\gamma)$ along the (unique) shortest
geodesic.
\par
Now the assertion immediately follows from the continuity on $U$ of the first
variation of the length in the direction of the field
$X(\gamma_*)=T_{\gamma}^{\gamma^*}(v(\gamma))$. This continuity 
follows from
the first variation formula for the length functional.
One just needs to perform easy calculations verifying
this continuity for all cases of elementary
mergers. That is,
it is necessary to consider the particular
cases, when $\gamma_i\longrightarrow \gamma$, where all $\gamma_i$
are of the same type, which can be obtained from (the type of) $\gamma$
by either 1) a splitting of a multiple point into two multiple points
(connected or not connected by a very short geodesic segment converging to
the point)
or into the pair ``multiple point-double point" (on one of the
segments adjacent to the multiple point);
or 2) a splitting of a double point into two double points connected
by a very short geodesic (converging to the double point). The general case 
of this formula follows by induction. We will omit the details of this 
straightforward verification.
\par {\bf 3.7.}
After these preliminaries we are going to prove that:
\par\noindent
A. There exists a sufficiently small positive $\tau_*\leq x$ such that
$g_{k,N}^{\tau_*}$ can be deformed to $\Gamma_k^0$; and
\par\noindent
B. For any positive $\tau\leq x$ there exists a deformation of $g_{k,N}^x$ to
$g_{k,N}^\tau$ inside $G_{k,N}^x$.
\par
In fact, as it was already noted, we are going to construct the deformation
by first definining a vector flow on an open subset of $(M^n)^{(N+1)k}$,
that includes the image of $G_{N,k}^x$ under the embedding discussed above.
This flow will be the same for both A and B. So the division
of our deformation into these two parts is somewhat artificial. Yet
we encounter different difficulties in these two situations:
When the length
is small, our main problem will be the lack of compactness at zero
length of the space $G_{N,k}^x\setminus G_{N,k}^0$, and we need to
prove that our flow decreases the length with a speed bounded from zero
by a constant. When the length is large, we do not want the distances
between double points that must be connected by a geodesic segment, to become
large, so from time to time we stop and perform the Birkhoff deformation.
\par {\bf 3.8.}
In order to prove A
we would first like
to establish a positive lower bound for $\Vert v(\gamma)\Vert$
for all $\gamma\in G_{k,N}^{\tau_0}\setminus G_{k,N}^0$ for a sufficiently
small positive $\tau_0$.
The key idea is to observe that this statement will be true for the
Euclidean space $R^n$ instead of $M^n$:
Assume that there exists a sequence of 1-cycles $\gamma_i$ in $R^n$ made of
at most $Nk$
straight line segments such that $\Vert v(\gamma_i)\Vert
\longrightarrow 0$.
Rescale $\gamma_i$ in $R^n$ 
so that  the maximal
length of an edge equals to one. (This does not affect $\Vert
v(\gamma_i)\Vert$.) Choose a convergent subsequence. Then its limit
must be a non-trivial stationary 1-cycle in $R^n$ of length $\geq 1$.
(Here we must check what happens with $\Vert v\Vert$ when an edge
collapses to a point
in the limit. It is easy to see that
$\Vert v(\lim_{i\longrightarrow\infty}(\gamma_i)\Vert\leq
2\lim\sup_{i\longrightarrow\infty} \Vert v(\gamma_i)\Vert$ 
in the situation, when we have a sequence of $\gamma_i$ of the same type,
and exactly one segment of $\gamma_i$ collapses to a point in the
limit. The number of such collapses is bounded from
above by $Nk-1$. Therefore the norm of the deformation vector of
the limit 1-cycle will be zero.)
But it is very easy to see that there are
no stationary 1-cycles in $R^n$.  So, we obtain a contradiction thereby
proving
the existence of a uniform positive lower bound for $\Vert v(\gamma)\Vert$
for all parametrized 1-cycles
that consists of at most $Nk$ straight line segments in $R^n$.
\par
If $\tau_0=\tau_0(M^n,N,k)$ is sufficiently small,
then any parametrized 1-cycle from
$G_{k,N}^{\tau_0}$
splits into several connected components contained in very small
balls in $M^n$. Applying the inverse of the exponential map we obtain
``almost" a 1-cycle in the tangent space to $M^n$ with ``almost" the
same angles. Now the existence of a uniform positive lower bound for
the norm of the deformation vectors of elements of $G_{k,N}$ for $R^n$
implies the existence of such uniform lower bound for all 1-cycles
from $G^{\tau_0}_{k,N}$.
\par
{\bf 3.9. Construction of a deformation of $g_{k,N}^{\tau_0}$ into
$\Gamma_k^0$.}
Now it is easy to find a countable set
$\{\gamma_l\}\subset G_{k,N}^{\tau_0}\setminus\Gamma_k^0$,
a locally finite covering of $G_{k,N}^{\tau_0}\setminus\Gamma_k^0$ by 
open balls $U_l$ centered at $\gamma_l$ and a subordinate partition of
unity that can be used to obtain a continuous function $\phi$ assigning
to every element $\gamma$ of
$G_{k,N}^{\tau_0}\setminus \Gamma_k^0$
a system of tangent vectors to $M^n$ at each
of its multiple or double points such that the variation of length of the cycles
in direction of $\phi(\gamma)$ is bounded from below by a positive
constant $\delta$. (In order to obtain a formal proof of the last assertion
we just need to establish the continuity of the first variation of the
length. This will be done in 3.10.)
\par
We construct this $\gamma_l$ and $U_l$ inductively
with respect to the partial order on types of elements of $G_{k,N}^x$
introduced in section 3.3.A above.
We start from elements of $G_{k,N}^x\setminus G_{k,N}^0$
of the highest possible type that form strata of a high codimension that are
closed in $G_{k,N}^x\setminus G_{k,N}^0$.
Construct a locally
finite covering (in $G_{k,N}^x$)
of the union of these strata, so that all centers
of open balls forming the covering are on the considered strata.
Automatically an open neighborhood of the union of these strata will be
covered.
Then we proceed to strata corresponding to types of the second highest
order. Points of closure of the union of
these strata in $G_{k,N}^x\setminus G_{k,N}^0$
that are not in these strata are in strata corresponding to the highest
type, and were already covered. Therefore we
can choose the covering so that no open ball from this covering
intersects an open neighborhood of strata corresponding to the highest
types that were previously covered. Also the centers of all balls should
be on the strata that are being covered.
We continue inductively in this way until we cover the whole
$G_{k,N}^x\setminus G_{k,N}^0$.
On each step of the inductive procedure
we  consider the union strata corresponding to the maximal types
that were not yet covered.
A neighborhood of the union of all
strata corresponding to higher types was already
covered on previous steps of the inductive procedure.
So we complete the covering of of the union of strata
that are being considered on the current step
by adding open balls centered at the considered strata
that have empty intersections with a neighborhood of the
union of all already covered
strata corresponding to all higher types.
\par
After the covering is completed, and a subordinate partition
of unity is chosen, we define
$\phi(\gamma)$ as the weighted sum of
$T_{\gamma_i}^{\gamma}(v(\gamma_i))$
over the set of indices $i$ such that $\gamma\in U_i$.
(Recall that components of $T_{\gamma_i}^\gamma(v(\gamma_i))$ are obtained
by the parallel transport of the corresponding components of $v(\gamma_i)$
along the shortest geodesics between corresponding multiple or double points
of $\gamma_i$ and $\gamma$. Of course, $U_i$ should be sufficiently small
in order for this definition to be unambigious, as we explained above.)
The weights are equal to the corresponding functions from the partition
of unity.
This construction
provides us with the flow $\Phi_t$ that deforms $g_{k,N}^\tau$ to $\Gamma_k^0$
in a finite time for each sufficiently small $\tau$ (as we will see below).
\par {\bf 3.10. The type cannot change during this deformation; the first
variation of the length is continuous.}
Let $p_1$, $p_2\in\gamma\in G_{k,N}^x\setminus G_{k,N}^0$ be
either two multiple points, or a multiple point
and a double point connected by a segment, or two double points 
connected by a segment.
Note that if the distance between $p_1$,
$p_2$
of $\gamma$ is very small, then $\gamma$ is in a small neighborhood of
a stratum corresponding to a higher type of cycles, where $p_1$
and $p_2$ merge into one point $p$. If this neighborhood is sufficiently
small, then {\it all} balls of the covering that cover $\gamma$ are
centered at strata corresponding to higher types, where $p_1$ and $p_2$
are merged into one point $p$. But then $p_1$ and $p_2$ will be deformed
(for some period of time) using the same vector fields on $M^n$.
This will be happening all the time while they will be
sufficiently close to each other. But since different integral trajectories
of a smooth vector field do not intersect, $p_1$ and $p_2$ cannot merge at least
until the moment of the deformation, when the total length of the 1-cycle
becomes zero.
\par
Also, observe that in the considered situation (when
$\gamma$ is close to a stratum corresponding to a higher type)
the first variation
of the length in the direction
of the vector field that we constructed
will be equal to a linear combination of variations in the direction
of vector fields of the form $T_{\gamma_l}^{\gamma}(v(\gamma_l))$,
considered in 3.6. (Here
$\gamma_l$ are located in strata corresponding to higher types, where
$p_1$ and $p_2$ merge into one point. 
The coefficients in the linear combination will be
the corresponding functions from the
partition of unity.) Therefore we can prove the continuity of the
first variation of length with respect to $\gamma$ as in 3.6. 
\par{\bf 3.11}
Note that for any $\gamma\in G_{k,N}^\tau$ $\Phi_t(\gamma)$ is defined only
until the moment of time $t(\gamma)$, when the length of $\Phi_t(\gamma)$
will become zero. But one can 
extend the domain of definition of $\Phi_t$ by defining
$\Phi_t(\gamma)=\Phi_{t(\gamma)}(\gamma)$ for $t>t(\gamma)$.
More precisely, we
take $\tau_*=\min\{x,\tau_0,\delta,inj(M^n)/4\}$,
where $\tau_0$ is as in 3.8, and $\delta$ is the lower bound
of the speed of decrease of the length introduced
at the beginning of 3.9, and just follow the flow
until we hit $\Gamma_k^0$. It is clear that:
1) For any element $\gamma\in g_{k,N}^\tau$ we will reach $\Gamma_k^0$ in time
 $t(\gamma)$ $\leq 1$; 
2)
Since the total length of $\gamma$ decreases,
the distance between any two points
on $M^n$ that should be connected by the shortest geodesic in order to
obtain $\Phi_t(\gamma)$ does not exceed $inj(M^n)/4$. Therefore $\Phi_t(\gamma)$
is unambigiously defined.
(Recall that we move multiple and double points
of $\gamma$ along trajectories of vector fields determined by the corresponding
components of $\phi(\gamma)$. In order to obtain $\Phi_t(\gamma)$ we connect
these points by the shortest geodesics.)
3) $t(\gamma)$ depends on $\gamma$ continuously (by virtue
of the implicit function theorem. In order to apply the
implicit function theorem we need to establish the continuous
differentiability of the length as a function of $\gamma$,
but it is equivalent to the continuity of the first variation of the
length of $\gamma$ in the direction of the vector field
that was established in the previous subsection.)
\par
Therefore the map assigning to $\gamma$
the point
$\Phi_{t(\gamma)}(\gamma)
\in \Gamma_k^0$, where the trajectory of the
flow reaches $\Gamma_k^0$ is continuous, and is the deformation of $G_{k,N}^\tau$ to $\Gamma_k^0$.
\par {\bf 3.12.}
It remains to prove the existence of a deformation of $g_{k,N}^x$ to
$g_{k,N}^{\tau_*}$ inside $\Gamma_k^x$.
Use the compactness of the closure $S$
of $G_{k,N}^x\setminus G_{k,N}^{\tau_*/2}$ to find
a finite open covering of $S$ by open neighborhoods $U_l$ of
$\gamma_l\in S$ such that for any $\gamma_*\in U_l$ the first variation
of the length of $\gamma$ in the direction of
$T_{\gamma_l}^{\gamma_*}(v(\gamma_l))$ does not exceed
$-\Vert v(\gamma_l)\Vert^2/2$. (Recall that we have proved the existence
of an open neighborhood $U$ with this property for any
$\gamma\in G_{k,N}^x\setminus\Gamma_k^0$.) Using a subordinate partition
of unity $\alpha_l(\gamma)$ define a vector field $\phi(\gamma)$ by
the formula
$\phi(\gamma)=\Sigma_l \alpha_l(\gamma)T_{\gamma_l}^{\gamma}(v(\gamma_l))$,
where we perform the summation only over indices $l$ such that $\gamma\in U_l$.
\par
Here the construction of this open covering is similar to that
in section 3.9: We construct it inductively starting from the strata
that correspond to the highest type and then proceed to cover strata
corresponding to lower types. On each step we add only open balls
centered at points on the considered strata that do {\it not} intersect
with an open neighborhood of the already covered strata corresponding
to higher types. As the result the $G_{k,N}^x$ type
will not be changing during the
considered stage of the deformation, and the proof of this fact 
coincides with the proof in section 3.10 above almost verbatim.
\par
Rescale $\phi(\gamma)$ by a continuous function equal
to zero on $G_{k,N}^{\tau_*/2}$ and to one on $G_{k,N}^{\tau_*}$. Denote
the resulting vector field by $\psi(\gamma)$.
Let $t_*=inj(M^n)/(16k)$. Consider the flow $\Psi_t(\gamma)$ defined
for all $\gamma\in g_{k,N}^x$ at $t\in [0,t_*]$ and determined by the vector
field $\psi(\gamma)$.
Our choice of $t_*$
guarantees that $\Psi_t(\gamma)$ will be in $G_{k,N}^x$.
(In other words, the distance between any pair of points
that need to be connected by a geodesic segment will not exceed $inj(M^n)/2$.) 
Observe that for any $\gamma\in g_{k,N}^x\setminus G_{k,N}^{\tau_*}$
the difference between the
length of $\gamma$ and the length of $\Psi_{t_*}(\gamma)$ will be at least
 $\delta t_*$,
where $\delta={1\over 2}\min_{i=1}^l\Vert v(\gamma_l)\Vert^2$.
\par
Now recall that the Birkhoff curve shortening process provides us with the
deformation $B_N$ of $\Gamma_k^x$ into $g_{k,N}^x$. The restriction of $B_N$
to $G_{k,N}^{x}\subset \Gamma_k^x$ is a deformation of
$G_{k,N}^{x}$ into
$g_{k,N}^{x}$ inside $\Gamma_k^x$. Apply $B_N$. The composition
of $B_N$ and $\Psi_{t_*}$ is a curve-shortening deformation of $g_{k,N}^x$ into
$g_{k,N}^{\max\{x-\delta t_*,\tau_*\}}$ inside $\Gamma_k^x$.
\par
Now we can apply $\Psi_{t_*}$ and then $B_N$ again and again, etc.
Let $K=[{1\over \delta t_*}]+1$. 
It is easy to see that $(B_N\Psi_{t_*})^K$ is the required
deformation of $g_{k,N}^x$ to $g_{k,N}^{\tau_*}$ inside $\Gamma_k^x$.
QED.
\par
{\bf Observation.} We would like to mention again (for a further use)
that the defined in 3.3 type of elements of $\Gamma_k^x$ does not
change during the considered deformation until possibly the very last
moment, when the length becomes zero.
\par
{\bf Remark.} It seems that one can find a more straightforward proof of
a modification of Proposition 4 suitable for our purposes using the same
approach as in the proof of Theorem 4.3 in [P] 
(but we had not checked the details). The
idea is to use the space of varifolds, the compactness of closed balls in
the space of varifolds (that follows from the Banach-Alaoglu theorem), and
the continuity of the first variation of the mass of a varifold in the direction
of a fixed vectors field (with respect to the variable varifold). Then, as in
the proof of Lemma 3, one constructs a map from the space
of varifolds of the bounded mass to the space of $C^1$-smooth vector fields
such that the mass decreases fast, when we apply the corresponding flow.
One can use an appropriate
locally finite covering of the space of varifolds of bounded non-zero mass
and a subordinate partition of unity in order to construct this map to
the space of $C^1$-smooth vector fields. It seems that in this way one can avoid
minor combinatorial complications in our proof related to 
the combinatorics of the space of parametrized 1-cycles made of piecewise
geodesics. 
\par
Finally observe that if $M^n$ is diffeomorphic
to $S^2$ we can combine Lemma 3 in the case of $k=2$ with the
observation made in the last paragraph of section 2 to
obtain an elementary proof of the following assertion used 
in our paper [NR]. (This assertion first appeared in [ClCo].)
\medskip
\proclaim{Proposition 5}
Let $M$ be a Riemannian manifold diffeomorphic to $S^2$.
Assume that for some $x$ there exists
a non-contractible map $f:S^1\longrightarrow\Gamma_2^x$.
Then there exists a non-trivial closed geodesic of length
$\leq x$ on $M$.
\endproclaim
\par
{\bf Proof.}
Lemma 3 implies the
existence of a non-trivial strongly stationary 1-cycle in
$\Gamma_2^x$. But, as we noted at the end of section 2,
each strongly stationary 1-cycle on a two-dimensional
manifold made of two segments is either
a closed geodesic or a union of two 
closed geodesics. But our stongly stationary 1-cycle is
non-trivial. Therefore either it is a non-trivial closed geodesic
or it contains a non-trivial closed geodesic as a subset.
QED.
\medskip
\centerline{\bf 4. Almgren correspondence}
\medskip
Now we are going to explain the Almgren correspondence  
in the 1-dimensional
case for 1-cycles made of finitely many closed curves.
(This simplified version is all that we will need for our purposes 
- see [A] for the full story).
Assume that we are given a continuous
map $A$ of a disk $D^m$ or a sphere $S^{m-1}$, or
more generally, a compact polyhedron $\vert K\vert$ into $\Gamma_k$.
The Almgren correspondence assigns to $A$ a $(dim\vert K\vert+1)$-
dimensional singular chain on $M^n$ as follows:
We can regard $\Gamma_k$ as a subset of the topological space of
all maps of the disjoint union of $k$ copies of $[0,1]$ into $M^n$.
To any of these maps we can assign a continuous map of
$X=D^m\times\bigcup_{i=1}^k [0,1]_i$
(or $X=S^{m-1}\times\bigcup_{i=1}^k[0,1]_i$,
or $X=\vert K\vert\times\bigcup_{i=1}^k[0,1]_i$) into $M^n$ in the standard
way. Further, since elements of $\Gamma_k$ are parametrized 
cycles we can identify
points of $2k$ sets $D^m\times \{0\}$, $D^m\times \{1\}$ (or, correspondingly
$2k$ sets $S^{m-1}\times \{0\}$, $S^{m-1}\times \{1\}$, or
$\vert K\vert\times \{0\}$, $\vert K\vert\times \{ 1\}$)
that are mapped into the same points
of $M^n$. The resulting quotient $X_A$ of $X$ can be quite complicated.
For our purposes we will only need the situation when $X_A$ can be triangulated
with a finite number of simplices.
Moreover, we are going to make an even stronger assumption that will
always hold when we apply the Almgren correspondance:
We assume that the compact polyhedron $K$
has a simplicial subdivision such that for any open simplex $\sigma$
of any dimension of this subdivision
all parametrized 1-cycles $A(t),\ t\in\sigma$ have the same type.
We will call this assumption a {\it local triviality} of $A$.
Now we can triangulate
$X_A$ so that the obvious projection $X_A\longrightarrow K$ becomes
a simplicial map
and consider the corresponding singular chain in $M^n$.
The local triviality of $A$ makes the following assertions
eveident:
If $K=S^{m-1}$ then that the resulting singular $m$-chain will
be a singular cycle. Its homology class does not depend on the chosen
triangulation of $X_A$. If $A$ is a map of $S^{m-1}$ to $\Gamma_k$
obtained as the restriction of a map $B$ of $D^m$ to $\Gamma_k$
then for any triangulation of $X_B$ the boundary of the corresponding
singular $(m+1)$-chain in $M^n$ will be the singular chain obtained from
$X_A\subset X_B$ with the induced triangulation. Therefore the singular
$m$-cycle in $M^n$ assigned to $A:S^{m-1}\longrightarrow M^n$ will
represent $0$ in $H_m(M^n)$ if $A$ is contractible.
\par
Note that if $f$ is a map of $S^{i-1}$ into $\Gamma_k^x$ satisfying
the local triviality assumption,
and $H:S^{i-1}\times [0,1]\longrightarrow \Gamma_k^x$ is the homotopy between
$f$
and a map $g$ of $S^{i-1}$ into $\Gamma_k^0$ constructed as in the proof of
Lemma 3 (in the situation when there are no non-trivial strongly stationary parametrized 1-cycles in $\Gamma_k^x$),
then $H$ also satisfies this assumption. This assertion immediately
follows from the observation made right after the proof of Lemma 3. 
This observation will easily imply the local triviality of $A$
in situations 
when we will need to use the Almgren correspondence in the course of proving
Theorems 1 and 2 in the next section.
\par
Yet for the sake of completeness we are also going to sketch how to modify the
above version
of the Almgren correspondence in order to make it work in the general case,
when we do not even have the triangulability of $X_A$. (This construction
will not be used in the present paper.)
Recall that at the beginning
of the proof of Lemma 3 in the previous section we introduced
$N=N(M^n,x)=[4x/inj(M^n)]+1$,
the spaces $g_{k,N}^x$ and $G_{k,N}^x$
made of parametrized 1-cycles of length $\leq x$
formed by $k$ piecewise geodesics made of at most $N$ geodesic
segments of length $\leq inj(M^n)/4$ and $\leq inj(M^n)/2$,
correspondingly, parametrized proportionally
to the arclength. We also defined the curve-shortening
Birkhoff deformation of $\Gamma_k^x$ into $g_{k,N}^x$.
Further recall that
$g_{k,N}^x$ and $G_{k,N}^x$ can be regarded as subsets
of $M^{kN}$. It is easy
to prove that there exists a subset $\bar g_{k,N}^x$
of $M^{kN}$ containing $g_{k,N}^x$
and contained in $G_{k,N}^x$ that
can be triangulated.
(The shortest formal way to prove the last assertion
is the following: Approximate the Riemannian metric on $M^n$
in $C^3$-topology by
an analytic Riemannian metric so that the distances on 
the resulting Riemannian manifold $\bar M^n$ do not exceed
corresponding distances on $M^n$. 
Now observe that metric balls of
radius $\leq inj(M^n)/2$ on $\bar M^n$ are subanalytic sets, the restriction
of the distance function to such metric
balls is a subanalytic function, and that
according to a well-known theorem of H. Hironaka subanalytic sets are
triangulable (cf. [B] for more details. See also [BM] for
the definition and basic properties of subanalytic sets and
function, including the proof of the mentioned theorem of
H. Hironaka.)
Therefore we can define
$\bar g_{k,N}^x$ as $g_{k,N}^x$ but using the distance function
on $\bar M^n$ instead of the distance function on $M^n$.)
It is easy to see that one can triangulate $\bar g_{k,N}^x$
so that the type of parametrized 1-cycles is constant on
every simplex of the triangulation.
Therefore we can take $x=\max_{y\in K} l(A(y))$,
compose $B_N$ with $A$, and take a simplicial approximation
$\bar A$ of the resulting composition
$B_N A: K\longrightarrow \bar g_{k,N}^x$.
Now we can consider the quotient $X_{\bar A}$ defined as above.
It is easy to see how to triangulate $X_{\bar A}$.
Therefore we can proceed as above,
and consider the corresponding singular chain in $M^n$.
It is clear that
if $K=S^{m-1}$, then this chain is a cycle, and this cycle represents $0\in
H_{m-1}(M^n)$ if and only if $A$ is contractible.  A small technical
complication that arises here is the following: Assume that we
apply the Almgren correspondence to a map $B:D^m\longrightarrow \Gamma_k^x$
and to the restriction $A$ of $B$ to $S^{m-1}=\partial D^m$. Then
the values of $x$ defined for these two mappings will , in general, be
different. Therefore the cycle in $M^n$ corresponding to
$A$ will not, in general, coincide with the boundary of
the chain corresponding to $B$. Yet, it is easy to contruct a homology between
these two cycles.
\par
In the case of a map of a
polyhedron into $Z_{(k)}$ one can consider a sufficiently fine subdivision
of the polyhedron $K$.
For any simplex of this subdivision the restriction of our
map onto this simplex lifts to $\Gamma_k$. Then we can proceed as in the
parametrized case (for this simplex). Finally sum the resulting chains over
all simplices of the triangulation. 
\par
As it was noted, F. Almgren  described in [A] a similar but somewhat
more technically complicated
construction applicable to maps of a polyhedron into
$\bar Z_1(M^n,Z)$ or even $\bar Z_k(M^n,Z)$.
\par
\medskip
{\bf 5. Proofs of Theorems 1 and 2.}
\medskip
{\bf Proof of Theorem 1.} Assume that $\alpha(M^n)$ and, in particular,
$l(M^n)$ is greater than $\frac{(n+2)!d}{3}$.  Consider a 
map $f: S^q \longrightarrow M^n$ representing a non-zero
element of $\pi_q(M^n)$, where $S^q$ is the standard sphere with 
a fine triangulation.  
Let $[S^q]$ be the fundamental 
class of $S^q$. 
Since the map is non-contractible, $f_*([S^q]) \neq 0\in H_q(M^n)$.
Let $D^{q+1}$ be a disc that has $S^q$ as its boundary.
Triangulate $D^{q+1}$ as the cone over the triangulation of $S^q$ (introducing
one new 0-dimensional simplex at the centre of $D^{q+1}$).
We will try to construct a singular $(q+1)$-chain in $M^n$, such 
that $f_*[S^n]$ will be 
its boundary, which is clearly impossible and will
result in a desired contradiction.

We are going to proceed inductively assigning an $i$-dimensional
singular chain in $M^n$ to each $i$-dimensional simplex of $D^{q+1}$
on the $i$-th step. This assignment will be denoted by $F$. 
The boundary of the singular chain that
corresponds to an arbitrary simplex $\sigma_i$ will be equal to the 
signed sum of chains assigned to simplices of the boundary of $\sigma_i$.
These signs will be the same as the signs 
with which the corresponding simplices 
enter $\partial \sigma_i$. These singular $i$-chains will be obtained from 
$(i-1)$-dimensional discs in the space of 1-cycles, particularly in
$Z_{(k(i-1))}$ for some function $k(i)$. One will use the 
Almgren correspondence between discs and chains explained in 
Section 4.  In turn, these $(i-1)$-dimensional discs are 
obtained by contracting $(i-2)$-dimensional spheres in $Z_{(k(i-2))}$
that are constructed from $(i-2)$-dimensional discs in
$Z_{(k(i-2))}$ corresponding to simplices of $\partial \sigma_i$ that
were constructed on the previous stage of induction.
\par
Alternatively, we can describe this procedure in the following
equivalent way (with somewhat more details):
We start from
a collection of maps of $D^1\longrightarrow Z_{(k(1))}$, where
$k(1)=3$, described below.
Then, inductively, for each simplex $\sigma^i$ in the considered triangulation
of $D^{q+1}$, ($i=2,3,\ldots, q+1$), we do the following:
\par\noindent
1) Construct a map of $S^{i-2}$ into $Z_{(k(i-1))}$
using $(i+1)$ maps of $D^{i-2}$ into $Z_{(k(i-2))}$
corresponding to $i+1$ $(i-1)$-dimensional simplices
in the boundary of $\sigma^i$
and obtained on the previous step of induction;
\par\noindent
2) We use our assumption about non-existence of sufficiently short non-trivial
strongly stationary 1-cycles and Proposition 4
to obtain a map of $D^{i-1}$
into $Z_{(k(i-1))}$ contracting this map of $S^{i-2}$. This map of
$D^{i-1}$ will correspond to $\sigma^i$.
\par
All these maps from discs and spheres into $Z_{(k(i))}$, $i=1,2,\ldots ,q$, can be lifted to $\Gamma_{k(i)}$.
(In fact, these maps will be obtained as compositions
of maps into $\Gamma_{k(i)}$ and projections $\Gamma_{k(i)}\longrightarrow Z_{(k(i))}$).
So, after the completion of this induction process for every $(q+1)$-dimensional simplex of the considered triangulation of $D^{q+1}$
we obtain a map from $D^{q}$ to $Z_{(k(q))}$, that can be lifted to $\Gamma_{k(q)}$.
Then we will apply the Almgren correspondence to each of the resulting maps of $D^q$ into $\Gamma_{k(q)}$
and sum the resulting singular $(q+1)$-chains in $M^n$. The result will be the required
$(q+1)$-chain.
\par
We will
begin with the {\bf $0$-skeleton} of $D^{q+1}\setminus S^q$ that consists of the
point $p$, the center of the disc.  We will assign to $p$ a singular 
$0$-chain that corresponds to an arbitrary
point $\tilde{p} \in M^n$. 
Now we will proceed to  the {\bf $1$-skeleton}: we will assign to the  
$1$-simplices of the form $[v_i,p]$ the singular $1$-chains that
correspond to minimal geodesics in $M^n$ that connect 
$\tilde{p}$ and $\tilde{v}_i=f(v_i)$ Next, we consider the 
{\bf $2$-skeleton}: Let $\sigma^2=[v_i,v_j,p]$ be a $2$-simplex of
$D^{q+1}\setminus S^q$. 
Consider
its boundary $\partial \sigma^2$ and the corresponding singular $1$-chain
on $M^n$, which equals to $[\tilde{v}_j,\tilde{p}]-[\tilde{v}_i,\tilde{p}]
+ [\tilde{v}_i,\tilde{v}_j]$.  
This can be viewed as 
a curve of length $\leq 2d+\epsilon$.  By our assumption, there 
is no closed geodesics of length smaller than or equal to $2d+\epsilon$,
so there is a curve shortening homotopy that connects this curve with
a point. Therefore, we assign to this $2$-simplex a singular $2$-chain
consisting of one singular 2-simplex that corresponds the surface generated
by this homotopy.
 The ``extension'' to  the
{\bf $3$-skeleton} will be somewhat different.  Let 
$\sigma^3=[v_{i_0},v_{i_1},v_{i_2}, v_{i_3}]$ be a $3$-simplex of
$D^{q+1}\setminus S^q$.
We want to find a singular $3$-chain to assign to  this simplex. 
Consider $\partial \sigma^3$.  There is a singular $2$-chain assigned to
the boundary of this simplex, which can also  be viewed as 
a $2$-sphere in $M^n$ of a particular shape.  Namely, to each
of the faces of the boundary not in $S^q$
there was  assigned a  surface generated by a 
curve shortening homotopy. Without any loss of generality we can assume
that the chosen fine
triangulation  of $S^q$ and the map of $S^q$ into $M^n$ were
chosen so that any two-dimensional simplex of the triangulation $S^q$ is also
mapped into the surface obtained by contracting its boundary in $M^n$ by
a homotopy that does not increase the length. 
As we will see, this 2-sphere corresponds to a $1$-sphere
in the space $Z_1(M^n, Z)$ that passes through the subset that
consists of cycles that are composed of no more than 12 curves, and
the mass of each such cycle is bounded from above by $8d+4\epsilon$. 
(See figure 1 to understand how this 1-sphere is constructed.)
In order to describe this correspondence 
let $e_1=[\tilde{v}_{i_0},\tilde{v}_{i_1}], 
e_2=[\tilde{v}_{i_0},\tilde{v}_{i_2}], e_3=[\tilde{v}_{i_0}, \tilde{v}_{i_3}],
e_4=[\tilde{v}_{i_1},\tilde{v}_{i_2}], e_5=[\tilde{v}_{i_1},\tilde{v}_{i_3}],
e_6=[\tilde{v}_{i_2},\tilde{v}_{i_3}]$, where each 
$[\tilde{v}_{i_s},\tilde{v}_{i_t}]$ is a minimal geodesic segment on the 
manifold.  Then we will let $\gamma_1=e_1+e_5-e_3, \gamma_2=-e_1+e_2-e_4,
\gamma_3=-e_2+e_3-e_6, \gamma_4=e_6-e_5+e_4$. Let $x_i$ be a point to 
which $\gamma_i$ contracts for $i=1,...,4$. Then the $1$-sphere in the
space of $1$-cycles will be constructed as follows:
let $\tilde{f}_i:D^2 \longrightarrow M^n, i=1,...,4$ be each of the 
four discs that make the $2$-sphere in $M^n$. Those discs correspond
to four maps $f_i:[0,1] \longrightarrow Z_1(M^n,Z)$, such that
$f_i(0)=T_{\{x_i\}}=0, f_i(1)=T_{\gamma_i}$.
These maps are precisely curve-shortening
homotopies used to obtain $\tilde f_i$; for any $t\in [0,1]$ $f_i(t)$
is a 1-cycle that consists of one closed curve. It can be regarded as
an element of $Z_{(3)}$ if we represent $\gamma_i$ as the collection
of three curves (=three sides of the triangle) glued at their endpoints,
and will keep track of these three curves during homotopies contracting
$\gamma_i$.
Now we will let $G_1:[0,1] \longrightarrow
Z_{(12)}$
be the map that for each $q \in [0,1]$ assigns $\Sigma_{i=1}^4 T_{f_i(q)}$,
(see figure 1(b)). 
Note that $G_1(0)=\Sigma_{i=0}^4 T_{\{x_i\}}$, which is the zero cycle,
(see figure 1(a)) and that $G_1(1)=\Sigma_{i=0}^4 T_{\gamma_i}$, which
is also the zero cycle, (see figure 1(c)). Thus, we obtain a map from 
$S^1$ to $Z_1(M^n, Z)$.

Proposition 4 implies that one of the following
is true about this $1$-sphere: either it can be contracted to a point 
without the mass increase, or there exists a stationary $1$-cycle of order $12$
of mass bounded by $8d + 4\epsilon$.  The existence of such a cycle for
all sufficiently small $\epsilon$ is precluded by our assumption.
So, there a disc that passes through $1$-cycles of order $12$ of mass
bounded by $8d+4\epsilon$. In order to apply Proposition 4 here we first
must check
that our map of $S^1$ into $Z^{8d+4\epsilon}_{(12)}$
is homotopic to a map that lifts to $\Gamma_{12}$.
The lifting of the map of $[0,1)\in S^1=[0,1]/{0,1}$ to $Z_{(12)}$ is obvious
, but we also need to find a homotopy of $\bigcup_{i=1}^4\gamma_i$ to
$\bigcup_{i=1}^4\{ x_i\}$, where each $\{ x_i\}$ is counted three times and is
regarded as a constant segment, in $\Gamma_{12}$.
This can be achieved by
first cancelling in a continuous way six pairs of edges $e_i$ with the
opposite orientations to a point, which is obviously possible (each pair
is connected over itself to the point corresponding to $t={1\over 2}$ counted
twice),
and then connecting $12$-tuples of
these points regarded as an element of $\Gamma_{12}$ with the constant
1-cycle $\{ x_1,x_2,x_3,x_4\}$ regarded as the cycle from $\Gamma_{12}$ (each
point is counted three times) using twelve
continuous paths. These paths follow our
homotopies restricted to the points of edges of $\gamma_i$ corresponding
to $t=0.5$ in the chosen parametrization of these edges.
As the result we obtain a lifting to $\Gamma_{12}$ of a map that differs from
$f$ only by a reparametrization. \par
Further, the last assertion of Proposition 4 implies that 
the lifting of our map $S^1\longrightarrow Z_{(12)}^{8d+4\epsilon}$ to
$\Gamma_{12}$ is also contractible. Indeed, we just need to verify the
contractibility of $12$ maps of $S^1\longrightarrow M^n$.
Of course, this fact follows
from the simply-connectedness on $M^n$.
However, there is even a more
straightforward geometric reason for contractibility of these $12$ circles
in $M^n$: each of them is formed by the trajectory of a homotopy $\tilde f_i$
from $x_i$ to a point in the middle of a geodesic segment $e_j$ traversed
two times in the opposite directions. 
\par
Using the Almgren correspondence we see that this disc corresponds 
to a $3$-chain that we will denote $\tilde{C}_{v_{i_0},...,v_{i_3}}$ 
in $M^n$ that has $F(\partial \sigma^3)$ as its boundary. So,
we will assign
$\tilde{C}_{v_{i_0},...,v_{i_3}}$ to the simplex $\sigma^3$. Now suppose we 
want to extend to {\bf $4$-skeleton}.
Consider any $4$-simplex of $D^{q+1}\setminus\partial D^{q+1}$
$\sigma^4=[v_{i_0},...,v_{i_4}]$.  
The following $3$-dimensional cycle in $M^n$:
$C_{v_{i_0},...,v_{i_4}}=\sum_{j=0}^{4}(-1)^j \tilde{C}_{v_{i_0},...,
\hat{v}_{i_j},...,v_{i_4}}$ corresponds to the boundary $\partial \sigma^4$ of
this simplex.  We claim that there exists a corresponding
map of the $2$-disc to $Z_1(M^n, Z)$, 
that takes the boundary of this disc
to the zero cycle, and such that the resulting sphere passes 
through the subset that consists of cycles of order $60$, i.e. cycles
that are composed of no more than twenty closed curves (each closed curve
being composed of three segments), of total
mass no more than $20 (2d+\epsilon)$. So, this map represents an element
of $\pi_2(Z_1(M^n))$, and this element will be precisely the element 
corresponding to the homology class of the considered 3-dimensional cycle
in $M^n$ under the Almgren isomorphism.
The above
map, denoted $G_2: \bar{D}^2 \longrightarrow Z_1(M^n, Z)$ 
will be constructed
as follows: let $f_j: \bar{D}^2 \longrightarrow Z_1(M^n, Z)$ be a map
corresponding to 
$(-1)^j\tilde{C}_{v_{i_0},...,\hat{v}_{i_j},...,v_{i_4}}$. Then
let $G_2(q)=\Sigma_{j=0}^4 T_{f_j(q)}$ for any $q \in \bar{D}^2$.

Let us now
examine $G_2(\partial \bar{D}^2)$.  We will see that for any $q \in 
\partial \bar{D}^2, G_2(q)$ will correspond to the union 
of $10$ pairs of closed curves, where
each pair will contain the same curve with two different orientations. In
other words, the corresponding 1-cycles will have opposite signs, and will
cancel. As the result we obtain a zero cycle.

Thus, we obtained
a $2$-sphere in the space of $1$-cycles.
We would like to apply Proposition 4. 
In order to do that we need to lift that map $G_2$ to the map
$\tilde{G}_2: \bar{D}^2 \longrightarrow \Gamma_{60}$ and 
examine what happens to the boundary of the disc under this map. Each 
point on the boundary is mapped to the $30$ pairs of segments in $M^n$.
Each pair consists of the same segment with opposite orientations.
We want to construct a homotopy between 
$\tilde{G}_2:\partial \bar{D}^2 \longrightarrow \Gamma_{60}$ and
a constant map, i.e. a map that will take a circle to $60$ point curves.
To construct this homotopy we 
cancel pairs of parametrized 1-cycles corresponding to
the 1-cycles with opposite orientations mentioned in
the previous paragraph in a continuous way.
We contract each pair $\gamma\bigcup -\gamma$
to $\gamma(0.5)$ over $\gamma$. Thus, we obtain a circle in the space
$\Gamma_{60}$, where each point $p \in S^1$ 
corresponds to $60$ constant paths (those paths are different for
$p \neq p'$). This circle can be interpreted as $60$ circles on $M^n$.
Since $M^n$ is simply connected, these circles can be contracted to an
arbitrary point in $M^n$.
After
contracting them we can obtain a point in the space $\Gamma_{60}$ 
made of $60$ constant segments.
So, combining $\tilde G_2$ with these two homotopies, we obtain a map
of the 2-disc into $\Gamma_{60}$ such that its boundary is mapped into
a point composed of $60$ constant segments. We can factor this
map through the sphere $S^2$ obtained from the disc by identifying
its boundary to a point (say, the north pole of the sphere.
In this case the southern hemisphere is mapped by $\tilde G_2$, and the
northern hemisphere is mapped into the subset of $\Gamma_{60}$ that corresponds
to the zero cycle (i.e. in $I^{-1}(0)$).)
So Proposition 4 applies: We have constructed a $2$-dimensional sphere
in the space $\Gamma_{60}$, and can conclude that
either this sphere can 
be contracted along the cycles of mass $\leq 20(2d+\epsilon)$, or
we have a stationary 1-cycle of mass controlled  from above by this bound.
(Here the verification of the contractibility of maps $F_j$ defined in the
text of Proposition 4 is equivalent to the contractibility of certain $60$
2-spheres in $M^n$. But now we are discussing the case of
$q\geq 3$, so $M^n$ is 2-connected.)
If $\epsilon$ is sufficiently small, then the second case is impossible.
In the former case, we obtain a $3$-disc in the space of $1$-cycles,
that corresponds to a $4$-chain that we will denote 
$\tilde{C}_{{v_0,...,v_4}}$. 
We will then assign this chain to $\sigma^4$. 

Now we can continue in the above manner until we fill the original
$q$-dimensional chain $f_*([S^q])$ by a $(q+1)$-dimensional chain in $M^n$.
As a corollary of our assumption nothing will stop us until we construct
the desired filling.
But as it was said before, this is impossible, and we obtain a contradiction
refuting our assumption. The constants $(q+2)!/3$ and $(q+2)!/2$ in the text of
Theorem 1 can be explained by the fact that all our 1-cycles
consist of at most  $4\times 5\times 6\times ... \times (q+2)$
closed curves
of length not exceeding $2d+\epsilon$, and
each of these closed curves consists of three segments. 
\par
Note that we can get a better estimate when $q=2$. In this case
we need to perform the extension process only till the dimension
$q+1=3$. We will need ``to represent" the union of four maps of $D^2$
to $M^n$ corresponding to four faces of a 3-dimensional simplex
as a map of a circle to $Z_1(M^n,Z)$. (Recall that these four maps
where obtained by contracting the maps of boundaries of these discs to
a point without increase of their lengths, see Fig. 1). In the body
of the proof we mapped a generic point $t\in [0,1]$ into the
1-cycle that corresponds to the union of four
curves obtained from homotopies contracting $\partial D^2_i$ at the moment $t$
(see Fig. 1(b)).
In the particular case $q=2$ we can proceed in a slighly different way. We can
start from two points obtained as the result of contraction of the maps of 
boundaries of $D^2_1$ and $D^2_2$ and to pass via cycles made of two
closed curves (obtained during the curve-shortening 
homotopies contracting the maps of $\partial D^2_1$
and $\partial D^2_2$) to the cycle made of the images of these two boundaries
(see Fig.1). The edge $[v_0,v_2]$ will be passed twice with opposite
orientation. Continue the homotopy by cancelling this edge. At the end of
this homotopy we obtain the map of the boundary of $D^2_1\bigcup D^2_2$.
Now note that $\partial (D^2_1\bigcup D^2_2) =\partial (D^2_3\bigcup D^2_4)$.
But we can similarly construct a homotopy between $\partial (D^2_3\bigcup
D^2_4)$ and the zero cycle that uses 1-cycles made of two curves obtained
from the curve-shortening homotopies contracting the maps of the boundaries
of $D^2_3$ and $D^2_4$. Joining these two homotopies we obtain the desired
homotopy between the zero 1-cycle and the zero 1-cycle, i.e. the desired circle
in the space of 1-cycles that passes through 1-cycles made of not more
than two closed curves of length not exceeding $2d+\epsilon$ (each).
See the proof of Theorem 1 in
[NR] for more details (in the situation when $M^n$ is diffeomorphic
to $S^2$. But this part of the proof is the same
there as in the more general situation.)
QED.
\medskip

{\bf Proof of Theorem 2.}
Assume $\alpha(M^n)$ and, in particular, $l(M^n)$ is greater than 
$2(n+1) n^n n!^{\frac{1}{2}}(n+2)! (vol M^n)^{\frac{1}{n}}$.
Then $\alpha(M^n)$ (and $l(M^n)$) are greater than $(n+2)!Fill Rad(M^n)$.
The definition of the filling radius implies that $M^n$ bounds in
the $(Fill Rad(M^n)+\delta)$-neighborhood of $M^n$ in $L^{\infty}(M^n)$.
Let $W$ ``fill" $M^n$ in the $(Fill Rad(M^n)+\delta)$-neighborhood of $M^n$
(that is $M^n=\partial W\ (mod\ 2)$: Since we do not assume that $M^n$ is
orientable, its fundamental homology class $[M^n]$ is defined only for
$Z_2$ coefficients.) 
Without any loss of generality we can assume that $W$ is a polyhedron.

Suppose $W$ together with $M^n$ is endowed with a very fine triangulation. 
We are going to try to 
construct a singular $(n+1)$-chain on $M^n$
such that the boundary of that chain is homologous to the boundary of $W$
(regarded as a chain).
 That is clearly impossible, so we will obtain a contradiction. 
 We will construct this chain by induction
with respect to the dimension of skeleta of $W$ .  
That is to each $i$-simplex of $W$ we will assing a singular 
$i$-chain on $M^n$.
We will begin with the {\bf $0$-skeleton} of $W$.  Let $v_i$ be a vertex of $W$.
Then $F(v_i)=\tilde{v}_i\in M^n=\partial W$, such that 
$d(v_i,\tilde{v}_i) =d (v_i, M^n) \leq Fill Rad M^n + \delta$.
Suppose $\tilde{v}_i, \tilde{v}_j$ come from the vertices 
$v_i, v_j$ of some simplex in $W$.  Then $d(\tilde{v}_i, \tilde{v}_j)
\leq 2 Fill Rad M^n + 3\delta$. (We assume here
that the triangulation of $W$ is fine so that the
lengths of 1-simplices of the triangulation are at most 
$\delta$.) Next, we are going to extend $F$
to the {\bf $1$-skeleton}.  We will assign to any 
$1$-simplex $[v_i,v_j] \subset W\setminus M^n$
a singular $1$-chain that corresponds 
to a minimal geodesic that connects $\tilde{v}_i$ and 
$\tilde{v}_j$ of length $\leq 2Fill Rad M^n +3\delta$.  Now we can see that 
the boundary of each $2$-simplex in $W$ is sent to a singular chain that
corresponds to a curve of length 
$\leq 6 Fill Rad M^n + 9\delta$, (we will assume that all simplices in $M^n$
are already short).

Next we are going to extend to the {\bf $2$-skeleton}.  Let $\sigma^2$ be
a $2$-simplex of $W$.  Consider its boundary $\partial \sigma^2$
and its corresponding singular $1$-chain.  There is a 
curve shortening homotopy that connects the curve corresponding to that chain
to a point.  So we will map
$\sigma^2$ to the chain that corresponds to the surface determined by this
homotopy. To ``extend"
$F$ to the {\bf $3$-skeleton} of $W$ consider an arbitrary 
$3$-simplex $\sigma^3$.  Consider  its boundary $\partial \sigma^3$ and the
corresponding singular $2$-chain, which can 
be viewed as $1$-sphere in the space $Z_1(M^n, Z)$ as in the proof of the 
Theorem 1. This sphere passes
through $1$-cycles of length $\leq 4 (6 Fill Rad M^n + 9 \delta)$.
Suppose this sphere cannot be contracted via the 1-cycles of smaller mass.
Then there exists minimal $1$-cycle of length 
$\leq 4 (6 Fill Rad M^n + 9\delta)$ contradicting our assumption.
(Here we use Proposition 4 from the previous section.
One can check that our spheres in the space of non-parametrized
1-cycles can be lifted
to spaces of parametrized 1-cycles exactly as this was done in the proof
of Theorem 1 above.)
So the above $1$-sphere can be
``filled'' by a disc that passes through $1$-cycles of mass 
not exceeding the above bound.  This disc corresponds to a singular $3$-chain
that has $F(\partial \sigma^3)$ as its boundary.  So we will assign this
chain to $\sigma^3$.
The procedure of ``extending" to $4$-skeleton
is similar  to the one in the proof of the Theorem 1:
At this point ``the image" of $\partial \sigma^4$ has been determined and 
it equals to 
$C_{v_{i_0},...,v_{i_4}} = (-1)^j
\sum_{j=0}^4 \tilde{C}_{v_{i_0},...,\hat{v}_{i_j},...,v_{i_4}}$.
This chain is in fact a 3-dimensional cycle in $M^n$,
and it can be interpreted as a sphere of dim $2$
in $Z_1(M^n, Z)$.
This sphere is constructed
as follows.  Let $f_j: \bar{D}^2 \longrightarrow Z_1(M^n)$ be a map
that corresponds to $(-1)^j \tilde{C}_{v_{i_0},...,\hat{v}_{i_j},...,v_{i_4}}$.
Then, let $G_2:\bar{D}^2 \longrightarrow Z_1(M^n)$ be a map that assigns
to every $q \in \bar{D}^2$ a point $\Sigma_{j=0}^4 T_{f_j(x)}$.  Then
it is easy to see that the boundary of the disc is mapped to the zero 1-cycle,
and we obtain a 2-sphere in $Z_{(60)}$.
Now we want to use
Proposition 4, so we need to lift this map of the 2-sphere to $\Gamma_{60}$.
First we lift $G_2$ in the obvious way and consider 
$\tilde{G}_2: \bar{D}^2 \longrightarrow \Gamma_{60}$. Next consider
what happens to $\tilde{G}_2: \partial \bar{D}^2 \longrightarrow \Gamma_{60}$.
We see that each point is mapped to the union of $30$ pairs of segments. Each
pair consists of the same segment with different orientation. Those segments
can be continuously contracted to their middles, (namely contract
$\gamma \cup -\gamma$ to $\gamma(0.5)$). Thus we obtain a homotopy
between the original circle and the circle that passes through constant
parametrized 1-cycles
only. This circle corresponds to $60$ circles on a manifold, which we want to
contract. 
Now unlike the proof of Theorem 1 above we cannot assume that $M^n$ is simply
connected. But we can contract these circles using the following simple
construction that will similarly work for all dimensions:
Consider the disc 
$\tilde{G}_2: \bar{D}^2 \longrightarrow \Gamma_{60}$. For any $p \in \bar{D}^2$
$\tilde{G}_2 (p)=\{ \gamma_1^p,...,\gamma_{60}^p \}$. For each
$p$ consider $\{ \gamma_1^p(0.5),...,\gamma_{60}^p(0.5) \}$.  This determines
a $2$-dimensional disc in $\Gamma_{60}$ that passes only through constant
parametrized 1-cycles.
This disc corresponds to $60$ discs on the manifold. Now our
circles can all be contracted over these discs, which establishes 
a homotopy between $\partial \bar{D}^2$ and a point in $\Gamma_{60}$.
Now we can construct the required map of the
$2$-dimensional sphere into $\Gamma_{60}$
providing the desired lifting to $\Gamma_{60}$ exactly as this was
done in the proof of Theorem 1. Also, note that in our present
situation maps $F_j$ defined in
the text of Proposition 4 are $60$ maps of $S^2$ to $M^n$ 
defined as follows: For $p$ in the southern hemisphere $F_j(p)$ is
defined as $(\tilde G_2(p))_j(0.5)$. (Here we identify the southern
hemisphere with $D^2$.)  For $p$ in the northern hemisphere between the
equator and a certain parallel $F_j(p)$ is constant on every meridian.
(This stage corresponds to contracting oppositely directed pairs
of segments to the points in the middle corresponding to $t=0.5$.)
Finally, $F_j$ maps the part of the northern hemisphere north of this
parallel using $(\tilde G_2(p))_j(0.5)$ again. So, up to
a homotopy $F_j$ maps both hemispheres of $S^2$ in the same way.
Hence $F_j$ is contractible by the obvious homotopy.
Therefore we can apply Proposition 4:
Since our 2-sphere passes only through sufficiently short 
parametrized 1-cycles
 and because of our assumption, it can be contracted through sufficiently
short parametrized 1-cycles. 
The 3-disc contracting this sphere corresponds to a 4-chain in 
$M^n$ filling the
3-cycle we started from. (Here we have used the Almgren correspondence 
explained
in the previous section). 
And so on.

It becomes obvious, that we can go on like that until we
``extend" $F$ to the $(n+1)$-skeleton of $W$ thereby obtaining
a $(n+1)$-singular chain in $M^n$ filling modulo $2$ the fundamental homology
class of $M^n$ which is clearly impossible. The resulting contradiction
proves the theorem. QED.
\par
{\bf Important Remark.} As the reader noticed, we had not really used
any notions or theorems of the geometric measure theory in order to
prove Theorems 1 and 2. In fact, we
could restrict ourselves to the consideration of  
rather elementary spaces $\Gamma_k$ of parametrized 1-cycles
instead of $Z_1(M^n,Z)$,
and to only use  Lemma 3 (but not Proposition 4),
since all maps to the space of (non-parametrized)
cycles that we constructed in the course of proving Theorems 1, 2  lift
to appropriate spaces $\Gamma_k$ (possibly after a reparametrization).
Yet we prefered to make our exposition as above in order to illustrate
geometric measure-theoretic origin of many ideas used in the present paper.

{\bf Acknowledgements.} The paper was partially written during the visit of 
the authors to IHES in Summer of 2001.
A part of this paper was written when both authors were working at
the Courant Insitute of Mathematical Sciences in 2001-2002.
The authors would like to thank
IHES  and the Courant Institute for their kind hospitality.
Both authors gratefully acknowledge the
partial support of their research by NSERC research grants. Regina Rotman
gratefully acknowledges the partial support of her research by the NSF
Post-doctoral Fellowship.

\newpage

\centerline{\epsfbox{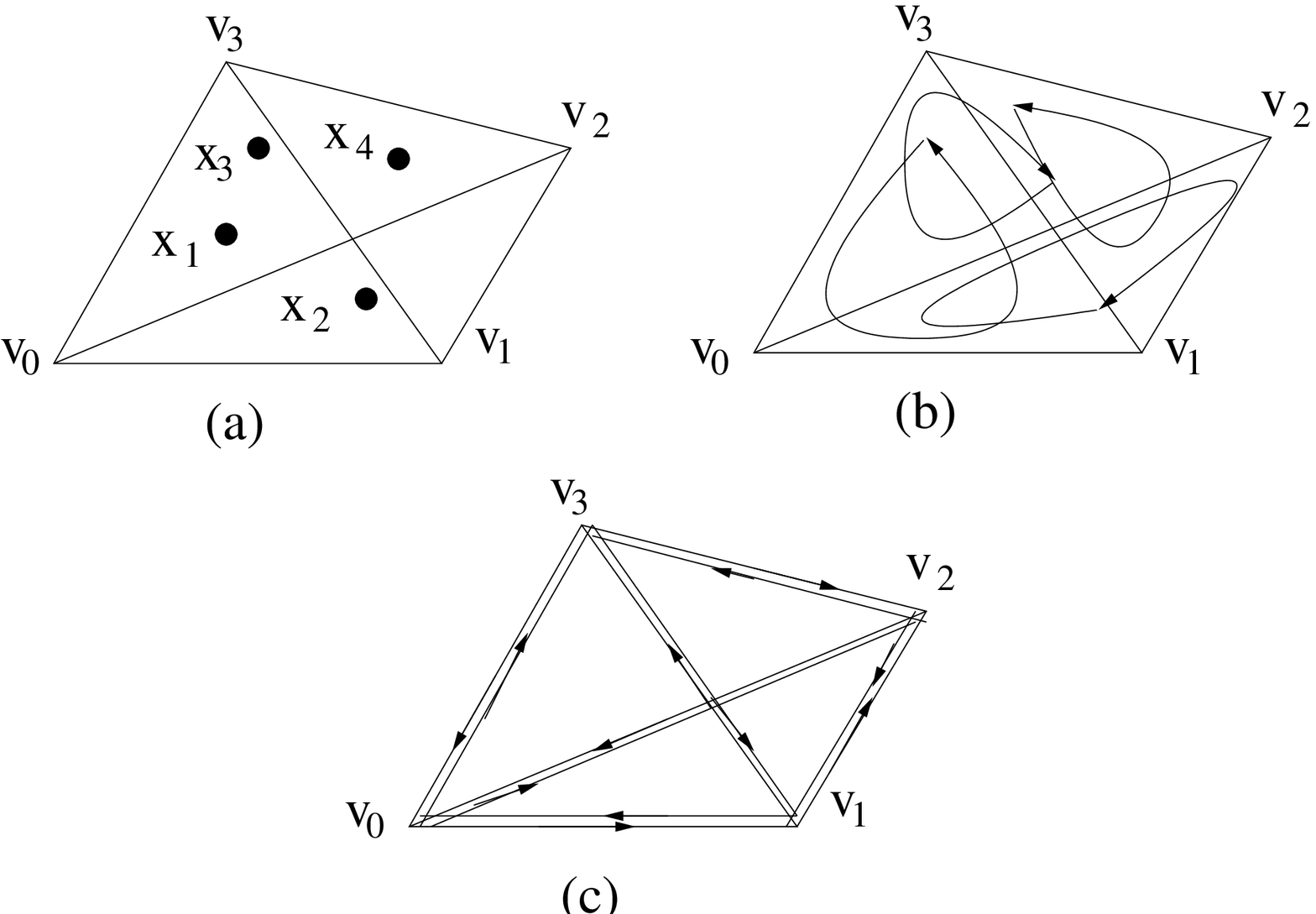}}

\enddocument